\documentclass[oneside]{amsart} 

\newcounter{version}
\newcounter{short}
\newcounter{long}
\usepackage{a4wide}
\usepackage{amsmath,amsthm,amscd}
\usepackage{colonequals}
\usepackage{enumerate}
\usepackage{amssymb,amsfonts} %
\input xy
\xyoption{all} \numberwithin{equation}{section}
\theoremstyle{definition}
\newtheorem{dfn}{Definition}[section]

\newtheorem{rem}[dfn]{Remark}
\theoremstyle{plain}
\newtheorem{thm}[dfn]{Theorem}
\newtheorem{prp}[dfn]{Proposition}
\newtheorem{cor}[dfn]{Corollary}
\newtheorem{lem}[dfn]{Lemma}

 \font\smallit=cmti8
\font\smallrm=cmr8

\newcommand\tens{{2.1}}
\newcommand\relvir{{5.1}}
\newcommand\comm{{5.2}}
\newcommand\vla{{6.1}}
\newcommand\Y{{6.2}}
\newcommand\heisrel{{6.3}}
\newcommand\glbrak{{6.4}}
\newcommand\virfa{{6.5}}
\newcommand\virfb{{6.6}}
\newcommand\glvir{{6.7}}
\newcommand\virfc{{6.8}}
\newcommand\actone{{6.9}}
\newcommand\acttwo{{6.10}}
\newcommand\actthree{{6.11}}
\newcommand\actfour{{6.12}}
\newcommand\bilpair{{8.1}}
\newcommand\copair{{8.2}}
\newcommand\bigpair{{8.3}}
\newcommand\contra{{8.4}}
\newcommand\contrb{{8.5}}
\newcommand\contrc{{8.6}}
\newcommand\contrd{{8.7}}
\newcommand\betavalue{{8.8}}
\newcommand\casim{{8.9}}
\newcommand\indu{{9.1}}
\newcommand\simpg{{9.2}}
\newcommand\Prz{{9.3}}
\newcommand\Eqq{{9.4}}
\newcommand\Pz{{9.5}}
\newcommand\gijk{{9.6}}
\newcommand\VW{{9.7}}
\newcommand\Pza{{9.8}}
\newcommand\Pzb{{9.9}}
\newcommand\oper{{9.10}}
\newcommand\oner{{9.11}}
\newcommand\Rza{{9.12}}
\newcommand\Rzb{{9.13}}
\newcommand\Rzc{{9.14}}
\newcommand\Rzd{{9.15}}
\newcommand\ghtwo{{9.16}}
\newcommand\sumtwo{{9.17}}
\newcommand\raq{{9.18}}
\newcommand\skews{{9.19}}
\newcommand\skewtwo{{9.20}}
\newcommand\skewtri{{9.21}}
\newcommand\zera{{9.22}}
\newcommand\zerb{{9.23}}
\newcommand\zerd{{9.24}}
\newcommand\zere{{9.25}}
\newcommand\zerf{{9.26}}
\newcommand\zerk{{9.27}}
\newcommand\zerh{{9.28}}
\newcommand\zerl{{9.29}}
\newcommand\zeri{{9.30}}
\newcommand\zerj{{9.31}}
\newcommand\zerm{{9.32}}
\newcommand\chone{{10.1}}
\newcommand\chtwo{{10.2}}
\newcommand\chthree{{10.3}}

\newcommand\C{{\mathbb C}}
\newcommand\Z{{\mathbb Z}}
\newcommand\Q{{\mathbb Q}}
\newcommand\T{{\mathbb T}}
\newcommand\vac{1 \hskip -1.0mm {\bf l}}

\newcommand\CC{{\mathcal C}}
\newcommand\D{{\mathcal D}}

\newcommand\K{\mathcal K}

\newcommand\PP{{\mathcal P}}
\newcommand\tP{{\widetilde{\PP}}}

\newcommand\g{{\frak g}}
\newcommand\dg{{\dot{\g}}}
\newcommand\Hyp{{\rm Hyp}}
\newcommand\Der{{\rm Der\,}}
\newcommand\gl{gl}
\newcommand\sln{sl}
\newcommand\dr{{\bf d}}
\newcommand\gln{\gl_N}
\newcommand\tsl{\widetilde{\raise 5pt \hbox{\vrule width 6pt height 1pt depth -2pt}} \kern-7pt {sl}}
\newcommand\stsl{\widetilde{\raise 3pt \hbox{\vrule width 6pt height 1pt depth -2pt}} \kern-7pt {sl}}
\newcommand\hgl{{\widehat{gl}}}
\newcommand\hsl{{\widehat{\sln}}}

\newcommand\hH{{\widehat{\mathcal H}}}
\newcommand\Vhyp{{V_{\hbox{\smallit Hyp}}^+}}
\newcommand\Vgl{{V_{\gl_N}}}
\newcommand\Lgl{{L_{\gl_N}}}
\newcommand\omgl{{\omega^{\gl_N}}}
\newcommand\omhei{{\omega^{\hbox{\smallit Hei}}}}
\newcommand\omvir{{\omega^{\hbox{\smallit Vir}}}}
\newcommand\omsl{{\omega^{\sln_N}}}
\newcommand\omhyp{{\omega^{{\hbox{\smallit Hyp}}}}}
\newcommand\om{{\omega}}
\newcommand\vir{{\hbox{\smallit Vir}}}
\newcommand\fer{{\hbox{\smallit fer}}}
\newcommand\bos{{\hbox{\smallit bos}}}

\newcommand\chara{{\hbox{\it char\,}}}

\newcommand\Vvir{{V_{\hbox{\smallit Vir}}}}
\newcommand\Lvir{{L_{\hbox{\smallit Vir}}}}
\newcommand\Cvir{{C_{\hbox{\smallit Vir}}}}
\newcommand\Mhyp{{M_{\hbox{\smallit Hyp}}}}
\newcommand\Mgl{{M_{\gl_N}}}
\newcommand\Mvir{{M_{\hbox{\smallit Vir}}}}

\newcommand\VZ{V_{{\mathbb Z}^N}}
\newcommand\hyp{{\hbox{\smallit Hyp}}}
\newcommand\Tr{{Tr}}
\newcommand\Cl{{Cl}}

\newcommand\len{\hbox{\rm len}}
\newcommand\sdeg{\hbox{\smallrm deg}}
\newcommand\pprime{{\prime\prime}}
\newcommand\ppprime{{\prime\prime\prime}}
\newcommand\gb{{\overline{g}}}
\newcommand\gh{{f}}
\newcommand\tw{{\widetilde{w}}}
\newcommand\tg{{\widetilde{g}}}
\newcommand\dd{\hbox{\bf d}}
\newcommand\vh{{\hbox{\rm v}_h}}
\newcommand\Id{{\hbox{\rm Id}}}
\newcommand\End{{\hbox{\rm End}}}

\newcommand\sgl{{\gl}}
\newcommand\ssl{{\sln}}
\newcommand\shei{{\hbox{\smallit Hei}}}
\newcommand\ot{{\otimes}}
\newcommand\darrow{\longrightarrow {\mkern -27mu} {\raise 6pt \hbox{\bf d}} {\mkern 16mu}}
\newcommand\marrow{{\hbox to 30pt{\rightarrowfill}}}
\newcommand\vrt{{\vrule height15pt width1pt}}
\newcommand\svrt{\raise5pt\hbox{\vrule height10pt width1pt}}

\newcommand\Ind{{\rm Ind}}
\newcommand\Map{{\rm Map}}

\newcommand\Venthree{{\rm Vect}{\T}^1}
\newcommand\VenN{{\rm Vect}{\T}^N}
\newcommand\Venone{{\rm Vect}{\T}^{N+1}}
\newcommand\Ventwo{{\rm Vect}{\T^2}}

\date{}
\begin{document}

\hskip 8.5cm
{\it To Yuri Alexandrovich Bahturin}

\

\title
[Representations of  Lie algebra of vector fields on a torus]
{Representations of Lie algebra of vector fields on a torus
and chiral de Rham complex}
\author{Yuly Billig}
\address{School of Mathematics and Statistics, Carleton University, Ottawa, Canada}
\email{billig@math.carleton.ca}
\author{Vyacheslav Futorny}
\address{ Instituto de Matem\'atica e Estat\'\i stica,
Universidade de S\~ao Paulo,  S\~ao Paulo,
 Brasil}
 \email{futorny@ime.usp.br}

\begin{abstract}
The goal of this paper is to study the representation theory  of a
classical infinite-dimensional Lie algebra -- the Lie algebra
$\VenN$ of vector fields on an $N$-dimensional torus for $N > 1$.
The case $N=1$ gives a famous Virasoro algebra (or its centerless
version - the Witt algebra).
 The algebra $\VenN$ has an important
class of tensor  modules parametrized by  finite-dimensional
modules of $\gl_N$. Tensor modules can be used in turn to
construct bounded irreducible modules for $\Venone$ (induced from
$\VenN$), which are the central objects of our study. We solve two
problems regarding these bounded modules: we  construct
their free field realizations and determine their characters. 
To solve these problems we analyze the structure of the
irreducible $\Omega^1 \left( \T^{N+1} \right) / d \Omega^0 \left(
\T^{N+1} \right) \rtimes \Venone$-modules constructed in \cite{B}. 
These modules remain
irreducible when restricted to the subalgebra $\Venone$, unless
they belongs to the {\it chiral de Rham complex}, introduced by
Malikov-Schechtman-Vaintrob \cite{MSV}.

\end{abstract}

\maketitle

\section{Introduction.}

In this paper we study the representation theory of a classical
infinite-dimensional Lie algebra -- the Lie algebra $\VenN$ of
vector fields on a torus. This algebra has a class of
representations of a geometric nature -- tensor modules, since
vector fields act on tensor fields of any given type via Lie
derivative. Tensor modules are parametrized by finite-dimensional
representations of $\gl_N$, with the fiber of a tensor bundle
being a $\gl_N$-module.

Irreducible $\gl_N$-modules yield tensor modules that are
irreducible over $\VenN$, with exception of the modules of
differential $k$-forms. In the latter case, the $\gl_N$-module is
irreducible, yet the modules of $k$-forms are reducible, which
follows from the fact that the differential of the de Rham complex
is a homomorphism of $\VenN$-modules. In the present paper we give
a vertex algebra analogue of this result.

 In case of a circle, a conjecture of Kac, proved by Mathieu \cite{M1}, states that for the Lie algebra of vector
fields on a circle an irreducible weight module with
finite-dimensional weight spaces is either a tensor module or a
highest/lowest weight module. There is a generalization of this
conjecture to an arbitrary $N$ due to Eswara Rao \cite{E2}. The
analogues of the highest weight modules in this case are defined
using the technique introduced by Berman-Billig \cite{BB}. These
modules are bounded with respect to one of the variables. It
follows from a general result of \cite{BB} that irreducible bounded
modules for the Lie algebra of vector fields on a torus have
finite-dimensional weight spaces, however the method of \cite{BB}
yields no information on the dimensions of the weight spaces. This
is the question that we solve in the present paper -- we find
explicit realizations of the irreducible bounded modules, using
which the dimensions of the weight spaces may be readily
determined.

 A partial solution of this problem for the $2$-dimensional torus was given by Billig-Molev-Zhang \cite{BMZ} using
non-commutative differential equations in vertex algebras. The
algebra of vector fields on $\T^2$ contains the loop algebra
$\tsl_2 = \C[t_0, t_0^{-1}] \otimes \sln_2$. This subalgebra plays
an important role in representation theory of $\Ventwo$. According
to the results of \cite{BMZ}, some of the bounded modules for
$\Ventwo$ remain irreducible when restricted to the subalgebra
$\tsl_2$. Futorny classified in \cite{Fu} irreducible generalized
Verma modules for $\tsl_2$. Such generalized Verma modules admit
the action of the much larger algebra $\Ventwo$.

 This relation between representations of $\tsl_2$ and $\Ventwo$ suggests that for the Lie algebra of vector
fields on $\T^N$, an important role is played by its subalgebra
$\tsl_N$. An unexpected twist here is that it is not the
generalized Verma modules for $\tsl_N$ that admit the action of
$\VenN$ for $N>2$, but rather the generalized Wakimoto modules.
The generalized Wakimoto modules are $\tsl_N$-modules that have
the same character as the generalized Verma modules, but need not
to be isomorphic to them.

 The generalized Wakimoto modules for $\tsl_N$ that we use here were constructed in \cite{B} in the context of
the representation theory of toroidal Lie algebras, however their special properties with respect to the loop
subalgebra $\tsl_N$ were not previously recognized.

 Let us outline the result of \cite{B} that we use here. Since one of the variables plays a special role, it will be more
convenient to work with an $(N+1)$-dimensional torus. To construct a full toroidal algebra, one begins
with the algebra of $\dg$-valued functions on $\T^{N+1}$:
$$ \Map (\T^{N+1}, \dg) \cong \C [t_0^{\pm 1}, t_1^{\pm 1}, \ldots, t_N^{\pm 1} ] \otimes \dg ,$$
where $\dg$ is a finite-dimensional simple Lie algebra. Next we take the universal central extension of this
multiloop algebra, with the center realized as the quotient of 1-forms on the torus by differentials of functions \cite{Ka}:
$$ \K = \Omega^1 \left( \T^{N+1} \right) / d \Omega^0 \left( \T^{N+1} \right) .$$
Finally, one adds the Lie algebra of vector fields on the torus:
$$\left(  \C [t_0^{\pm 1}, \ldots, t_N^{\pm 1} ] \otimes \dg \oplus \K \right) \rtimes \Venone .$$
Irreducible bounded modules for this Lie algebra were constructed
in \cite{B} using vertex algebra methods. Note that the results of
\cite{B} admit a specialization to $\dg = (0)$. The multiloop
algebra then disappears, leaving behind, like the smile of the
Cheshire Cat, the space of its central extension:
$$ \K \rtimes \Venone .$$
It turns out that it is easier to study representations of this Lie algebra, rather than of vector fields alone, because
of the duality between vector fields and $1$-forms. Representation theory of this larger Lie algebra is controlled
by a tensor product of three vertex algebras:
$$\Vhyp \otimes \Vgl \otimes \Vvir , $$
a subalgebra of a hyperbolic lattice vertex algebra, the affine
$\hgl_N$ vertex algebra at level $1$ and the Virasoro vertex
algebra of rank $0$. The tensor product of the first two
components, $\Vhyp \otimes \Vgl$ is one of the bounded modules for
$\K \rtimes \Venone$, and in fact it is a generalized Wakimoto
module for the subalgebra $\tsl_{N+1} \subset  \Venone$. Then the
results of \cite{BMZ} suggest that there is a chance that $\K
\rtimes \Venone$-modules constructed in \cite{B} remain
irreducible when restricted to $\Venone$. The study of this
question is the main part of the present paper. The answer that we
get is remarkably parallel to the classical picture with the
tensor modules. We prove that a bounded irreducible  $\K \rtimes
\Venone$-module remains irreducible when restricted to the
subalgebra of vector fields, unless it belongs to the {\it chiral
de Rham complex}, introduced by Malikov-Schechtman-Vaintrob
\cite{MSV} (for arbitrary manifolds).

It is only in very special situations an irreducible module remains irreducible when restricted to a subalgebra.
A prime example of this is the basic module for an affine Kac-Moody algebra, which remains irreducible
when restricted to the principal Heisenberg subalgebra. This exceptional property of the basic module leads to
its vertex operator realization and is at heart of several spectacular applications of this theory.

The space of the chiral de Rham complex is the vertex superalgebra
$$\Vhyp \otimes \VZ ,$$
where $\VZ$ is the lattice vertex superalgebra associated with the standard euclidean lattice $\Z^N$.
The vertex superalgebra $\VZ$ is graded by fermionic degree:
$$\VZ = \mathop\oplus\limits_{k\in \Z} \VZ^k ,$$
and the components
$$\Vhyp \otimes \VZ^k $$
are irreducible $\K \rtimes \Venone$-modules. Yet for these
modules the restriction to the Lie algebra of vector fields is no
longer irreducible since the differential of the chiral de Rham
complex
$$\dd : \; \Vhyp \otimes \VZ^k \rightarrow \Vhyp \otimes \VZ^{k+1}$$
is a homomorphism of $\Venone$-modules.

 In fact it was noted in \cite{MS1} that these components admit the action of the Lie algebra
$\C [t_0, t_0^{-1}] \otimes \VenN$, but here we prove a much
stronger result.

 In conclusion, we make two curious observations. The Lie algebra of vector fields on a torus has a trivial
center, yet its representation theory is described in terms of vertex algebras $\Vhyp$ and $\Vgl$ that involve non-trivial central extensions. The central charges of these tensor factors cancel out to give a vertex algebra of total rank $0$.

The final remark is that the chiral de Rham complex is an essentially super object, whereas the Lie algebra
of vector fields we started with, is classical.

The structure of the paper is as follows. We introduce the main objects of our study 
in Sections 2, 3 and 4. We discuss vertex algebras and their applications to the representation theory of $\Venone$ in Sections 5 and 6. In Section 7 we introduce 
the generalized Wakimoto modules for the loop algebra $\tsl_N$. We construct 
a non-degenerate pairing for the bounded modules in Section 8. We prove the main result on irreducibility in Section 9 and make a connection with the chiral de Rham complex in the final section of the paper.

\

\

\section{Lie algebra of vector fields and its tensor modules}

 We begin with the algebra of Fourier polynomials on an $N$-dimensional torus $\T^N$.
 Introducing the variables $t_j = e^{ix_j}, \, j=1,\ldots,N$, we
realize the algebra of functions as Laurent polynomials
$\C[t_1^{\pm 1}, \ldots, t_N^{\pm 1}]$.
 The Lie algebra of vector fields on a torus is
 $$\VenN = \Der \C [t_1^{\pm 1}, \ldots, t_N^{\pm 1}] =
 \mathop\oplus\limits_{p=1}^N  \C [t_1^{\pm 1}, \ldots, t_N^{\pm 1}] {\partial \over \partial t_p} .$$
 It will be more convenient for us to work with the degree derivations
 $d_p = t_p {\partial \over \partial t_p}$ as the free generators of $\VenN$ as
 a $\C[t_1^{\pm 1}, \ldots, t_N^{\pm 1}]$-module:
 $$\VenN =
 \mathop\oplus\limits_{p=1}^N  \C [t_1^{\pm 1}, \ldots, t_N^{\pm 1}] d_p .$$

 The Lie bracket in $\VenN$ is then written as
 $$[t^r d_a, t^m d_b ] = m_a t^{r+m} d_b - r_b t^{r+m} d_a , \quad a,b = 1,\ldots, N.$$
 Here we are using the multi-index notations $t^r = t_1^{r_1} \ldots t_N^{r_N}$ for
 $r = (r_1,\ldots,r_N) \in \Z^N$.

 The Cartan subalgebra $\left< d_1, \ldots, d_N \right>$ acts on $\VenN$ diagonally and induces
 on it a $\Z^N$-grading.

 The Lie algebra of vector fields (on any manifold) has a class of representations of a geometric nature.
 Vector fields act via Lie derivative on the space of tensor fields of a given type. The resulting tensor modules are
 parametrized by representations of $\gl_N$. Let us describe the construction of tensor modules in case
 of a torus $\T^N$.

 Fix a finite-dimensional $\gl_N$-module $W$. 
In case when $W$ is irreducible, the identity matrix acts as multiplication by a scalar
$\alpha \in \C$.  Let $\gamma\in \C^N$. We define the tensor module
 $T = T(W,\gamma)$ to be the vector space
 $$ T = q^\gamma \C[q_1^{\pm 1},\ldots, q_N^{\pm 1}] \otimes W $$
 with the action given by
 $$t^r d_a (q^\mu \otimes w) = \mu_a q^{\mu+r} \otimes w + \sum_{p=1}^N r_p q^{\mu+r} \otimes E^{pa} w ,
\eqno{(\tens)}$$
where $r\in\Z^N, \mu\in \gamma + \Z^N, \, a=1, \ldots, N$ and $E^{pa}$ is the matrix with $1$ in $(p,a)$-position
and zeros elsewhere.

\begin{thm} \label{DRc} [\cite{E1}, cf. \cite{Ru}] Let $W$ be an irreducible
finite-dimensional $\gl_N$-module. The tensor module $T(W,\gamma)$
is an irreducible $\VenN$-module, unless it appears in the de Rham
complex of differential forms
$$q^{\gamma} \Omega^0(\T^N) \darrow q^{\gamma} \Omega^1(\T^N) \darrow \; \ldots \; \darrow q^{\gamma} \Omega^N(\T^N) .$$
The middle terms in this complex are reducible $\VenN$-modules,
while the terms \break $q^{\gamma} \Omega^0(\T^N)$ and $q^{\gamma}
\Omega^N(\T^N)$ are reducible whenever $\gamma\in\Z^N$.
\end{thm}

 Note that de Rham differential $\dr$ is a homomorphism of $\VenN$-modules.
Let us specify irreducible $\gl_N$-modules that correspond to the tensor modules in de Rham complex.
The modules of functions
 $\Omega^0$ and the module of differential $N$-forms $\Omega^N$ correspond to 1-dimensional $gl_N$-modules
 $W$ on which the identity matrix acts as multiplication by $\alpha = 0$ and $\alpha = N$ respectively.
 The remaining modules $\Omega^k, \, k=1,\ldots,N-1,$ are the highest weight modules for $\sln_N$ with the fundamental
 highest weights $\omega_k$ and $\alpha = k$ (see e.g. \cite{E1}). Even though
 they correspond to irreducible $\gl_N$-modules,  tensor modules of differential forms are reducible since the kernels and images of the differential $\dr$ are obviously the submodules in $\Omega^k$.

\

\section{Bounded modules}

Our goal is to generalize to an arbitrary $N$ the category of the
highest weight modules over $\VenN$. In our constructions one of
the coordinates will play a special role. From now on, we will be
working with the $N+1$-dimensional torus and will index our
coordinates as $t_0, t_1, \ldots, t_N$, where $t_0$ is the
``special variable''. We would like to construct modules for the
Lie algebra $\D = \Venone$ in which the ``energy operator'' $-d_0$
has spectrum bounded from below.

Let us consider a $\Z$-grading of $\D$ by degrees in $t_0$. This
$\Z$-grading induces a decomposition
$$ \D = \D_- \oplus \D_0 \oplus \D_+$$
into subalgebras of positive, zero and negative degrees in $t_0$.
The degree zero part is
$$\D_0 = \mathop\oplus\limits_{p=0}^N \C [t_1^{\pm 1}, \ldots, t_N^{\pm 1}] d_p .$$
In particular, $\D_0$ is a semi-direct product of the Lie algebra of
vector fields on $\T^N$ with an abelian ideal $\C [t_1^{\pm 1}, \ldots,
t_N^{\pm 1}] d_0$.

We begin the construction of a bounded module by taking a tensor
module for $\D_0$. Fix a finite-dimensional irreducible
$\gl_N$-module $W$, $\beta \in \C$ and $\gamma\in\C^N$. We define a
$\D_0$-module $T$ as a space
$$T = q^\gamma \C [q_1^{\pm 1}, \ldots, q_N^{\pm 1}] \otimes W$$ with the tensor
module action (\tens) of the subalgebra $\VenN \subset \D_0$ and
with \break $\C [t_1^{\pm 1}, \ldots, t_N^{\pm 1}] d_0$ acting by
shifts
$$ t^r d_0 (q^\mu \otimes w) = \beta \, q^{\mu+r} \otimes w .$$
Next we let $\D_+$ act on $T$ trivially and define $M(T)$ as the
induced module
$$M(T) = \Ind_{\D_0 \oplus \D_+}^{\D} T \cong U(\D_-) \otimes T.$$
The module $M(T)$ has a weight decomposition with respect to the
Cartan subalgebra $\left< d_0, \ldots, d_N \right>$ and the (real
part of) spectrum of $-d_0$ on $M(T)$ is bounded from below. However
the weight spaces of $M(T)$ that lie below $T$ are all
infinite-dimensional.

It turns out that the situation improves dramatically when we pass
to the irreducible quotient of $M(T)$. One can immediately see
that the Lie algebra $\D$ belongs to the class of Lie algebras
with polynomial multiplication (as defined in \cite{BB}), whereas
tensor modules belong to the class of modules with polynomial
action. A general theorem of \cite{BB} (see also \cite{BZ}) yields
in this particular situation the following

\begin{thm} \label{TBB} (\cite{BB}) (i) The module $M(T)$ has a unique maximal
submodule $M^{rad}$.

(ii) The irreducible quotient $L(T) = M(T)/M^{rad}$ has
finite-dimensional weight spaces.

\end{thm}

This leads to the following natural questions:\\
\

{\bf Problem 1.} Determine the character of $L(T)$.

{\bf Problem 2.} Find a realization of $L(T)$.

\

In \cite{BMZ} these problems were solved for some of the modules
$L(T)$ in case of a 2-dimensional torus ($N=1$).  In the present
paper we will give a solution in full generality for any $N$.

 \

 \section{Toroidal Lie algebras}

 For a finite-dimensional simple Lie algebra $\dg$ we consider a
 multiloop algebra
\break
$\C[t_0^{\pm 1}, \ldots, t_N^{\pm 1}] \otimes \dg$.
Its universal central extension has a realization with center
 $\K$ identified as the quotient space of 1-forms by differentials
 of functions \cite{Ka},
 $$\K = \Omega^1 (\T^{N+1}) / d\Omega^0 (\T^{N+1}) .$$
The Lie bracket in
$$ \C[t_0^{\pm 1}, \ldots, t_N^{\pm 1}] \otimes \dg \oplus \K$$
is given by
$$ [f_1(t) \otimes g_1, f_2(t) \otimes g_2] =
f_1 f_2 \otimes [g_1, g_2] + (g_1 | g_2) \overline{ f_2 d f_1 }, $$
where $g_1, g_2 \in \dg$,
$f_1, f_2 \in  \C[t_0^{\pm 1}, \ldots, t_N^{\pm 1}]$,
$(\cdot | \cdot)$ is the Killing form on $\dg$ and
$\overline{\raise5pt\hbox{\hskip 0.35cm}}$ denotes the projection
$\Omega^1 \rightarrow \Omega^1 /d \Omega^0$.

 We will set 1-forms $k_a = t_a^{-1} dt_a$, $a=0, \ldots, N$ as
generators of $\Omega^1(\T^{N+1})$ as a free $\C[t_0^{\pm 1},
\ldots, t_N^{\pm 1}]$-module. We will use the same notations for
their images in $\K$.

 The Lie algebra $\D = \Venone$ acts on the universal central extension
 of the multiloop algebra with the natural action on
$ \C[t_0^{\pm 1}, \ldots, t_N^{\pm 1}] \otimes \dg$, and the action
on $\K$ induced from the Lie derivative action of vector fields on
$\Omega^1$:
$$ f_1 (t) d_a (f_2(t) k_b) = f_1 d_a (f_2) k_b + \delta_{ab} f_2
d(f_1), \; \; a,b = 0, \ldots, N.$$
The full toroidal Lie algebra is a semi-direct product
$$\g = \left( \C[t_0^{\pm 1}, \ldots, t_N^{\pm 1}] \otimes \dg \oplus \K
\right) \rtimes \D .$$
In fact \cite{B} treats a more general family of Lie algebras, where the
Lie bracket in $\g$ is twisted with a 2-cocycle
$\tau\in H^2 (\D, \Omega^1 /d\Omega^0)$. However for the purposes of the
present work we need to consider only the semi-direct product, i.e., set
$\tau = 0$.

 A category of bounded modules for the full toroidal Lie algebra is studied
in \cite{B} and realizations of irreducible modules in this
category are given.
 The constructions of \cite{B} admit a specialization $\dg = (0)$, which yields
representations of the semi-direct product
$$\D \ltimes \K.$$
The approach of the present paper is to look at the
representations of this semidirect product, constructed in
\cite{B}, and to study their reductions to the subalgebra $\D$ of
vector fields on $\T^{N+1}$. Surprisingly, as we shall see below,
most of the irreducible modules for $\D \ltimes \K$ remain
irreducible when restricted to $\D$.

 In order to describe here the results of \cite{B}, we will need to present a
background material on vertex algebras.

\

\section{Vertex superalgebras: definitions and notations}

Let us recall the basic notions of the theory of the vertex
operator (super) algebras. Here we are following \cite{K} and
\cite{L}.

\begin{dfn}  A vertex superalgebra is a $\Z_2$-graded vector space
$V$ with a distinguished vector $\vac$ (vacuum vector) in $V$, a
parity-preserving
operator $D$ (infinitesimal translation) on the space $V$, and a
parity-preserving linear map $Y$ (state-field correspondence)
$$
Y(\cdot,z): \quad V \rightarrow (\End V)[[z,z^{-1}]],$$ $$  a
\mapsto Y(a,z) = \sum\limits_{n\in\Z} a_{(n)} z^{-n-1} \quad
(\hbox{\rm where \ } a_{(n)} \in \End V),$$

 such that the
following axioms hold:

\noindent
(V1) For any $a,b\in V, \quad a_{(n)} b = 0 $ for $n$ sufficiently large;

\noindent
(V2) $[D, Y(a,z)] = Y(D(a), z) = {d \over dz} Y(a,z)$ for any $a \in V$;

\noindent
(V3) $Y(\vac,z) = \Id_V \, z^0$;

\noindent
(V4) $Y(a,z) \vac \in V [[z]]$ and $Y(a,z)\vac |_{z=0} = a$ for any $a \in V$
\ (self-replication);

\noindent
(V5) For any $a, b \in V$, the fields $Y(a,z)$ and $Y(b,z)$ are mutually local, that is,
$$ (z-w)^n \left[ Y(a,z), Y(b,w) \right] = 0, \quad \hbox{\rm for \ }
 n  \hbox{\rm \ sufficiently large} .$$

A vertex superalgebra $V$ is called a vertex operator superalgebra (VOA) if, in addition,
$V$ contains a vector $\omega$ (Virasoro element) such that

\noindent
(V6) The components $L_n = \omega_{(n+1)}$ of the field
$$ Y(\omega,z) = \sum\limits_{n\in\Z} \omega_{(n)} z^{-n-1}
= \sum\limits_{n\in\Z} L_n z^{-n-2} $$
satisfy the Virasoro algebra relations:
$$  [ L_n , L_m ] = (n-m) L_{n+m} + \delta_{n,-m} {n^3 - n \over 12}
\Cvir,
\eqno{(\relvir)}$$
where $\Cvir$ acts on $V$ by scalar, called the {\it rank} of $V$.

\noindent
(V7) $D = L_{-1}$;

\noindent
(V8) Operator $L_0$ is diagonalizable on $V$.
\end{dfn}

This completes the definition of a VOA.

As a consequence of the axioms of the vertex superalgebra
we have the  following important commutator formula:
$$\left[ Y(a,z_1), Y(b,z_2) \right] =
\sum_{n \geq 0} {1 \over n!}  Y(a_{(n)} b, z_2) \left[ z_1^{-1}
\left( {\frac{\partial}{{\partial} z_2}} \right)^n \delta \left(
{\frac{z_2}{z_1}} \right) \right] . \eqno{(\comm)}$$
 As usual, the delta function is
$$ \delta(z) = \sum_{n\in\Z} z^n .$$
By (V1), the sum in the left hand side of the commutator formula
is actually finite. The commutator in the left hand side of (\comm)
is of course the supercommutator.

Let us recall the definition of a normally ordered product of two fields.
For a formal field $a(z) = \sum\limits_{j\in \Z} a_{(j)} z^{-j-1}$ define its positive and
negative parts as follows:
$$ a(z)_- = \sum_{j=0}^\infty a_{(j)} z^{-j-1}, \quad
a(z)_+ = \sum_{j=-1}^{-\infty} a_{(j)} z^{-j-1} .$$
Then the normally ordered product of two formal fields $a(z), b(z)$ of parities
$p(a), p(b) \in \{ 0, 1 \}$ respectively, is defined as
$$ :a(z) b(z): = a(z)_+ b(z) + (-1)^{p(a)p(b)} b(z) a(z)_- \; .$$

The following property of vertex superalgebras will be used extensively in this paper:
$$ Y(a_{(-1)} b, z) =\; : Y(a,z) Y(b,z): \;  , \quad \hbox{\rm for all \ } a,b\in V .$$

%

%


\

\section{Vertex Lie superalgebras}
 An important source of the vertex superalgebras is provided by the vertex Lie superalgebras.
In presenting this construction we will be following \cite{DLM}
(see also  \cite{FKRW}, \cite{K}, \cite{P}, \cite{Ro}).

Let $\mathcal L$ be a Lie superalgebra with the basis $\{ u_{(n)},
c_{(-1)} | u\in \mathcal U, c\in\CC, n\in\Z \}$ ($\mathcal U$,
$\CC$ are some index sets). Define the corresponding fields in
$\mathcal L [[z,z^{-1}]]$:
$$ u(z) = \sum_{n\in\Z} u_{(n)} z^{-n-1}, \quad c(z) = c_{(-1)} z^0, \quad
u\in\mathcal U, c\in\CC .$$ Let $\mathcal F$ be a subspace in
$\mathcal L [[z,z^{-1}]]$ spanned by all the fields $u(z), c(z)$
and their derivatives of all orders.

\begin{dfn}
 A Lie superalgebra $\mathcal L$ with the basis as above is called a vertex Lie superalgebra
if the following two conditions hold:

(VL1) for all $x, y \in \mathcal U$,
$$ [x(z_1), y(z_2) ] = \sum\limits_{j=0}^n f_j(z_2)
\left[ z_1^{-1} \left(\frac{\partial}{{\partial} z_2} \right)^j
\delta \left( {z_2 \over z_1} \right) \right], \eqno{(\vla)}$$
where $f_j(z) \in\mathcal F, n \geq 0$ and depend on $x, y$,

(VL2) for all $c\in\CC$, the elements $c_{(-1)}$ are central in
$\mathcal L$.
\end{dfn}

\


Let $\mathcal L_{(+)}$ be a subspace in $\mathcal L$ with the
basis $\{ u_{(n)} \big| u\in\mathcal U, n\geq 0 \}$ and let
$\mathcal L_{(-)}$ be a subspace with the basis $\{ u_{(n)},
c_{(-1)} \big| u\in\mathcal U, c\in\CC, n<0 \}$. Then $\mathcal L
= \mathcal L_{(+)} \oplus \mathcal L_{(-)}$ and $\mathcal L_{(+)},
\mathcal L_{(-)}$ are in fact subalgebras in $\mathcal L$.

The universal enveloping vertex algebra $V_{\mathcal L}$ of a
vertex Lie superalgebra $\mathcal L$ is defined as the induced
module
$$V_{\mathcal L} = \Ind_{\mathcal L_{(+)}}^\mathcal L (\C \vac) = U(\mathcal L_{(-)}) \ot \vac,$$
where $\C \vac$ is a trivial 1-dimensional $\mathcal L_{(+)}$
module.

\begin{thm} \label{voa} [\cite{DLM}, Theorem 4.8]  Let $\mathcal L$ be
a vertex Lie superalgebra. Then

(a) $V_{\mathcal L}$ has a structure of a vertex superalgebra with
the vacuum vector $\vac$, infinitesimal translation $D$ being a
natural extension of the derivation of $\mathcal L$ given by
$D(u_{(n)})$  $=$ $-n u_{(n-1)}$, $D(c_{(-1)}) = 0$, $u\in\mathcal
U$, $c\in\CC$, and the state-field correspondence map $Y$ defined
by the formula:
$$Y \left( a^1_{(-1-n_1)} \ldots a^{k-1}_{(-1-n_{k-1})} a^k_{(-1-n_k)} \vac, z
\right) $$
$$ = \;
:\left( {1\over n_1 !} \left( \frac{\partial}{{\partial} z}
\right)^{n_1} a^1 (z) \right) \ldots :\left( \frac{1}{n_{k-1} !}
\left( \frac{\partial}{{\partial} z} \right)^{n_{k-1}} a^{k-1} (z)
\right) \left( {1\over n_{k} !} \left( \frac{\partial}{{\partial}
z} \right)^{n_k} a^k (z) \right): \ldots : \quad , \eqno{(\Y)}$$
where $a^j \in \mathcal U, n_j \geq 0$ or $a^j \in\CC, n_j =0$.

(b) Any bounded $\mathcal L$-module is a vertex superalgebra
module for $V_{\mathcal L}$.

(c) For an arbitrary character $\chi: \CC \rightarrow \C$, the quotient module
$$ V_{\mathcal L} (\chi) =  U(\mathcal L_{(-)}) \vac /  U(\mathcal L_{(-)}) \big< (c_{(-1)} - \chi(c))
\vac \big>_{c\in\CC}$$
is a quotient vertex superalgebra.

(d) Any bounded $\mathcal L$-module in which $c_{(-1)}$ act as
$\chi(c) \Id$, for all $c\in\CC$, is a vertex superalgebra module
for $V_{\mathcal L}(\chi)$.

\end{thm}
The value $\chi(c)$ is referred to as {\it central charge} or {\it level}.

\

The vertex algebra that controls representation theory of $\D
\ltimes \K$ is the tensor product of three VOAs: a subalgebra
$\Vhyp$ of a hyperbolic lattice vertex algebra, an affine $\hgl_N$
vertex algebra $\Vgl$ at level 1, and the Virasoro vertex algebra
$\Vvir$ of rank 0. In order to apply the results of \cite{B} to
representation theory of $\D \ltimes \K$, we use specializations
$\dg = (0)$ and $\tau = 0$. In this specialization one has to fix
the following values for the various central charges appearing in
(\cite{B}, Theorems 5.3 and 6.4):
$$ c = 1, \; \; c_{sl_N} = 1; \; \; c_{\hbox{\smallit Hei}} = N, $$
$$ c_{\hbox{\smallit VH}} = {N \over 2}, \; \; c_{\hbox{\smallit Vir}}^\prime = 0.$$

Let us briefly review the constructions of these three vertex
algebras.

{\it Hyperbolic lattice VOA.}
 Consider a hyperbolic lattice $\Hyp$,
which is a free abelian group on $2N$ generators $\{ u^a, v^a | a =
1, \ldots, N \}$ with the symmetric bilinear form
$$ ( \cdot | \cdot ) : \; \Hyp \times \Hyp \rightarrow \Z, $$
defined by
$$ (u^a | v^b) = \delta_{ab}, \; \; (u^a | u^b) = (v^a | v^b) = 0,
\; \; a,b = 1, \ldots, N. $$
 We complexify $\Hyp$ to get a $2N$-dimensional vector space
 $$\mathcal H = \Hyp \otimes_{\mathbb Z} \C$$
 and extend the bilinear form by linearity on $\mathcal H$. Next, we
 affinize $\mathcal H$ by defining a Heisenberg Lie algebra
 $$\hH = \C [t, t^{-1}] \otimes \mathcal H \oplus \C K$$
 with the Lie bracket
 $$ [x_{(n)}, y_{(m)}] = n (x|y) \delta_{n,-m} K, \;\; x,y \in\mathcal H , \eqno{(\heisrel)}$$
 $$ [K, \hH] = 0 .$$
 Here and in what follows, we are using the notation $x_{(n)} = t^n \otimes x$.

 The algebra $\hH$ has a triangular decomposition
 $\hH = \hH_- \oplus \hH_0 \oplus \hH_+$, where $\hH_0 = 1\otimes \mathcal H \oplus \C K$,
 and $\hH_\pm = t^{\pm 1} \C[t^{\pm 1}] \otimes \mathcal H$.

 Let $\Hyp^+$ be an isotropic sublattice of $\Hyp$ generated by $\{
 u^a | a=1,\ldots, N \}$. We consider its group algebra $\C[\Hyp^+]
 = \C [e^{\pm u^1},\ldots,e^{\pm u^N}]$ and define the action of
$\hH_0 \oplus \hH_+$ on $\C[\Hyp^+]$ by
$$x_{(0)} e^y = (x|y) e^y, \;\; K e^y = e^y, \;\; \hH_+ e^y = 0.$$
To be consistent with our previous notations, we set $q_a = e^{u^a}$,
$a =1, \ldots, N$.

Finally, let $\Vhyp$ be the induced module
$$\Vhyp = \Ind_{\hH_0 \oplus \hH_+}^{\hH} \left( \C[\Hyp^+]\right) .$$
We coordinatize $\Vhyp$ as a Fock space over $\hH$:
$$\Vhyp =  \C [q_1^{\pm 1}, \ldots, q_N^{\pm 1}] \otimes
 \C [u_{pj}, v_{pj} |^{p = 1, \ldots, N}_{j = 1,2, \ldots}]  ,$$ where
$\hH$ acts by operators of multiplication and
differentiation:
$$ u^p_{(-j)} = j u_{pj}, \;\; u^p_{(j)} = {\partial \over \partial
v_{pj}}, \;\; u^p_{(0)} = 0, $$
$$ v^p_{(-j)} = j v_{pj}, \;\; v^p_{(j)} = {\partial \over \partial
u_{pj}}, \;\; v^p_{(0)} = q_p {\partial \over \partial q_p}, $$ for $p
= 1, \ldots, N, \; j = 1, 2, \ldots .$

The space $\Vhyp$ has the structure of a vertex algebra - it is a
vertex subalgebra in the vertex algebra corresponding to lattice
$\Hyp$. We give here the values of the state-field
correspondence map on the generators of this vertex algebra:
$$Y (u_{p1}, z) = u^p (z) = \sum_{j\in\Z} u^p_{(j)} z^{-j-1} ,$$
$$Y (v_{p1}, z) = v^p (z) = \sum_{j\in\Z} v^p_{(j)} z^{-j-1} , \;\; p = 1, \ldots, N,$$
$$Y (q^r, z) = q^r exp
\left( \sum_{p=1}^N r_p \sum_{j=1}^\infty {z^j \over j} u^p_{(-j)}
\right) exp \left( - \sum_{p=1}^N r_p \sum_{j=1}^\infty { z^{-j}
\over j} u^p_{(j)} \right) .$$ The Virasoro element of $\Vhyp$ is
$$\omega^\hyp = \sum_{p=1}^N u_{p1} v_{p1} $$
and the Virasoro field is
$$Y (\omega^\hyp, z) = \sum_{p=1}^N :u^{p}(z) v^{p}(z): \; . $$
The rank of $\Vhyp$ is $2N$.

Fix $\gamma\in\C^N$. The space
$$M_\hyp (\gamma) = q^\gamma \C [q_1^{\pm 1}, \ldots, q_N^{\pm 1}] \otimes
\C [u_{pj}, v_{pj} |^{p = 1, \ldots, N}_{j = 1,2, \ldots}]$$ has a
natural structure of a simple module for $\Vhyp$ (see e.g.
\cite{BB} for details).

\

{\it Affine $\hgl_N$ VOA.} The second vertex algebra that we will
need is the affine $\hgl_N$ vertex algebra at level 1. Since $\gl_N$
is reductive, but not simple, it has more than one affinization.
Here we consider a particular version of $\hgl_N$:
$$\hgl_N = \C[t, t^{-1}] \otimes \gl_N \oplus \C C$$
with the Lie bracket
$$ [t^n \otimes X, t^m \otimes Y] = t^{n+m} \otimes [X,Y] + n
\delta_{n, -m} \Tr (XY) C, \;\; X,Y\in\gl_N . \eqno{(\glbrak)}$$


 We note that $\hgl_N$ is a vertex Lie algebra and consider its universal
enveloping vertex algebra $\Vgl$ at level 1 (i.e., $\chi(C) = 1$).

Let us give the value of the
state-field correspondence map on the generators of this affine
vertex algebra:
$$ Y( X_{(-1)} \vac, z) = X(z) = \sum_{j\in\Z} X_{(j)} z^{-j-1}, \;\;
\hbox{ \rm \ for \ } X\in\gl_N .$$ Since $\gl_N$ has a decomposition
$\gl_N = \sln_N \oplus \C I$, where $I$ is the identity $N \times N$
matrix, the affine $\hgl_N$ vertex algebra is the tensor product of
the affine $\hsl_N$ vertex algebra and a Heisenberg vertex algebra.
The Virasoro element $\omgl$ of $\Vgl$ can be thus written as a sum
of the Virasoro elements $\omsl$ for the affine $\hsl_N$ vertex
algebra and $\omhei$ for the Heisenberg vertex algebra. The usual
formula for the Virasoro element in affine vertex algebra gives the
following explicit expression:
$$\omsl = {1 \over 2 (N+1)} \left( \sum_{i,j=1}^N E^{ij}_{(-1)} E^{ji}_{(-1)}
\vac - {1\over N} I_{(-1)} I_{(-1)} \vac \right) . \eqno{(\virfa)}$$
 The rank of the affine $\hsl_N$ vertex algebra at level 1 is $N-1$.

For the Heisenberg vertex algebra we choose a non-standard Virasoro
element (see \cite{B}, (4.33)):
$$ \omhei = {1\over 2N} I_{(-1)} I_{(-1)} \vac + {1\over 2} I_{(-2)} \vac .\eqno{(\virfb)}$$
The rank of this Heisenberg VOA is $1-3N$.

Adding the two Virasoro elements, we get the Virasoro element for
$\Vgl$:
$$\omgl = {1\over 2(N+1)} \left( \sum_{i,j=1}^N E^{ij}_{(-1)} E^{ji}_{(-1)} \vac
+ I_{(-1)} I_{(-1)} \vac \right) + {1\over 2} I_{(-2)} \vac . \eqno{(\glvir)}$$
The corresponding Virasoro field is
$$Y(\omgl, z) = {1\over 2(N+1)} \left( \sum_{i,j=1}^N :E^{ij} (z) E^{ji}
(z):  + :I(z) I(z): \right) + {1\over 2} {d \over dz} I(z) . \eqno{(\virfc)}$$ The
rank of $\Vgl$ is $-2N$.

Let $W$ be a finite-dimensional simple module for $\gl_N$.
Let $C$ act on $W$ as the identity operator and set
$\left( t \C[t] \otimes \gl_N \right) W = 0$.
Construct the generalized Verma module for the Lie algebra $\hgl_N$
as the induced module from $W$, and consider its irreducible quotient
$\Lgl (W)$.  Then $\Lgl (W)$ is a simple module for the vertex algebra
$\Vgl$.

\

{\sl Virasoro VOA.} The last vertex algebra that we need to
introduce is the Virasoro vertex algebra $\Vvir$ of rank 0. The
Virasoro Lie algebra (\relvir) is a vertex Lie algebra with
$\mathcal U = \{ \omvir \}$ and $\CC = \{ \Cvir \}$, where
$$\omvir (z) =  \sum_{j\in\Z} \omvir_{(j)} z^{-j-1} = \sum_{j\in\Z} L_j z^{-j-2} .$$
Let $\Vvir$ be its universal enveloping vertex algebra of zero central charge, $\chi(\Cvir) = 0$.


Let  $L_\vir (h)$ be the irreducible highest weight module for the Virasoro
Lie algebra with central charge $0$ with the highest weight vector $\vh$,
satisfying $L_0 \vh = h \vh$.

\

 The vertex algebra that controls representation theory of $\D
 \ltimes \K$ is the tensor product of the sub-VOA $\Vhyp$ of the
 hyperbolic lattice vertex algebra, affine $\hgl_N$ vertex
 algebra $\Vgl$ at level 1, and the Virasoro VOA $\Vvir$ of rank 0
 $$  \Vhyp \otimes \Vgl \otimes \Vvir $$
 with the Virasoro element
 $$\omega = \omhyp + \omgl + \omvir .$$
 The rank of this VOA is $2N - 2N + 0 = 0$. Now we are ready to present
 a result of \cite{B} (Theorems 5.3 and 6.4):

 \begin{thm} \label{thDK} (\cite{B})
  (i) Let $\Mhyp$, $\Mgl$, $\Mvir$ be modules for
 $\Vhyp$, $\Vgl$ and $\Vvir$ respectively. Then the tensor product
 $$\mathcal M = \Mhyp \otimes \Mgl \otimes \Mvir$$
 is a module for the Lie algebra $\D \ltimes \K$ with the action given as follows:
 $$ \sum_{j\in\Z} t_0^j t^r k_0 z^{-j} =
 k_0 (r,z) \mapsto Y(q^r, z), \eqno{(\actone)}$$
 $$ \sum_{j\in\Z} t_0^j t^r k_a z^{-j-1} =
k_a (r,z) \mapsto   u^a (z) Y(q^r, z), \eqno{(\acttwo)}$$
 $$ \sum_{j\in\Z} t_0^j t^r d_a z^{-j-1} =
 d_a (r,z) \mapsto : v^a (z) Y(q^r, z): + \sum_{p=1}^N
 r_p E^{pa}(z) Y(q^r, z), \eqno{(\actthree)}$$
 $$ \sum_{j\in\Z} t_0^j t^r d_0 z^{-j-2} =
 d_0 (r,z) \mapsto - :Y(\omega,z) Y(q^r,z): - \sum_{i,j=1}^N
 r_i u^j(z) E^{ij}(z) Y(q^r, z)$$
$$ {\hskip 4cm} + \sum_{p=1}^N r_p \left( {d \over dz} u^p(z) \right)
Y(q^r, z) , \eqno{(\actfour)}$$
for $a = 1, \ldots, N$.

(ii) The module
$$L(W,\gamma,h) = \Mhyp(\gamma) \otimes \Lgl (W) \otimes \Lvir(h)$$
is an irreducible module over the Lie algebra $\D \ltimes \K$.

\end{thm}

\section{Generalized Wakimoto modules.}

In \cite{BMZ} the structure of irreducible modules $L(T)$ over the
Lie algebra of vector fields was determined in case of a
2-dimensional torus ($N=1$). It turned out that the situation was
analogous to the case of a basic module for an affine Kac-Moody
algebra, which remains irreducible when restricted to the
principal Heisenberg subalgebra \cite{LW}, \cite{KKLW}. For the
Lie algebra of vector fields on $\T^2$ this role is played by its
loop subalgebra $\tsl_2 = \C[t_0, t_0^{-1}] \otimes \sln_2$.
Indeed we have $\sln_2$ embedded into $\Venthree$:
$$\sln_2 \cong \left< t_1^{-1} d_1, d_1, t_1 d_1 \right> \subset
\C[t_1, t_1^{-1}] d_1 .$$ This extends to an embedding
$$ \tsl_2 \cong \C[t_0, t_0^{-1}] \otimes
\left< t_1^{-1} d_1, d_1, t_1 d_1 \right> \subset \Ventwo .$$

The following theorem was proved in \cite{BMZ}:

\begin{thm} \label{BMZ} (\cite{BMZ}) Let $W$ be the 1-dimensional $\gl_1$-module in
which the identity matrix acts as multiplication by $\alpha \in
\C$. Assume $\alpha \not\in \Q$, $\beta = {\alpha (\alpha-1) \over
2}$,
 and let $\gamma \in \C$. Then the module $L(T) = L(\alpha, \beta,
\gamma)$ over the Lie algebra $\Ventwo$ remains irreducible when
restricted to subalgebra $\tsl_2$.
\end{thm}

{Note that \cite{BMZ} uses a different convention for the sign of
$\alpha$.}

The loop algebra $\tsl_2 = \C[t_0, t_0^{-1}] \otimes \sln_2$ is
$\Z$-graded by degree in $t_0$. This gives its decomposition $\tsl_2
= \tsl_2^+ \oplus \tsl_2^0 \oplus \tsl_2^-$. The zero part $\tsl_2^0
\cong \sln_2$ is a subalgebra in $\C [t_1, t_1^{-1}] d_1$ and thus
acts on $T$. The positive part $\tsl_2^+$ acts on $T$ trivially. We
can form the generalized Verma module over $\tsl_2$:
$$ \Ind_{\stsl_2^0 \oplus \stsl_2^+}^{\stsl_2}
 T(\alpha, \beta, \gamma)  \cong U(\tsl_2^-) \otimes
T(\alpha, \beta, \gamma).$$ By the results of \cite{Fu}, this
generalized Verma module is irreducible over $\tsl_2$ if and only
if $\alpha\not\in {1\over 2} \Z$. This gives the following

\begin{cor} \label{cBMZ} (\cite{BMZ}) Let $\alpha \in \C,
\alpha\not\in\Q$, $\beta= {\alpha (\alpha-1) \over 2}$. Then the
$\Ventwo$-module $L(\alpha, \beta, \gamma)$ when restricted to
$\tsl_2$ is isomorphic to the generalized Verma module over
$\tsl_2$:
$$L(\alpha, \beta, \gamma) \cong U(\tsl_2^-) \otimes T(\alpha, \beta,
\gamma) .$$
\end{cor}

 These results show that the loop subalgebra $\tsl_2$ plays a
 crucial role in representation theory of the Lie algebra of vector
 fields on $\T^2$. It is natural to conjecture that in the
 representation theory of $\D = \Venone$ such a role is played by
 the loop algebra $\tsl_{N+1}$. Indeed, $\D_0$ has $\VenN$ as
 a subalgebra and $\VenN$ contains $\sln_{N+1}$ (see e.g. \cite{M2}).
Thus $\Venone$ contains the loop algebra $\tsl_{N+1} =
 \C  [t_0, t_0^{-1}] \otimes \sln_{N+1}$. The modules $T(W, \beta, \gamma)$
 may be viewed as $\sln_{N+1}$-modules, and we can form the generalized
 Verma module over $\tsl_{N+1}$:
 $$ U(\tsl_{N+1}^-) \otimes T(W, \beta, \gamma) .$$
 It turns out, however, that in general, the action of $\tsl_{N+1}$
 on the generalized Verma module can not be extended to the action
 of the bigger algebra $\Venone$. Instead one should use certain
 generalized Wakimoto modules. We define the generalized Wakimoto
 modules in the following way:

 \begin{dfn} Let $T$ be an $\sln_{N+1}$-module. A {\it
 generalized Wakimoto module} $\mathcal M$ with top $T$ is an
 $\tsl_{N+1}$-module that contains $T$ as an $\sln_{N+1}$-submodule
 with $\tsl_{N+1}^+$ acting on $T$ trivially and having the property
 that the character of $\mathcal M$ coincides with the character of the
 generalized Verma module for $\tsl_{N+1}$:
 $$ \chara \mathcal M = \chara \left( U(\tsl_{N+1}^-) \right)
 \times \chara T .$$
\end{dfn}

 The generalized Verma module is by definition a generalized
 Wakimoto module. In case when the generalized Verma module is
 irreducible, it is the only generalized Wakimoto module with the
 given top $T$. As we mentioned above, for $\tsl_2$ this happens for
 the tops $T(\alpha, \beta, \gamma)$ with $\alpha\not\in
 {1\over 2}\Z$ \cite{Fu}. For $N >1$ and any top $T(W,\beta,\gamma)$ with
 a finite-dimensional $\gl_N$-module $W$, the generalized Verma module
 over $\tsl_{N+1}$ is always reducible.

 We shall now see that Theorem \ref{thDK} yields
 a  construction of a generalized Wakimoto module for $\tsl_{N+1}$.
 These modules admit the action of the whole algebra of vector fields on
 $\T^{N+1}$.

\begin{prp} \label{genWak} Let $M_{\gl_N} (W)$ be the generalized
Verma module for $\hgl_N$ at level $1$, induced from an
irreducible finite-dimensional ${\gl_N}$-module $W$. Then the
module
$$\mathcal M = q^\gamma \C [q_1^{\pm 1}, \ldots, q_N^{\pm 1}]  \otimes \C[u_{pj}, v_{pj}] \otimes M_{\gl_N} (W)$$
is a generalized Wakimoto module for the loop algebra $\tsl_{N+1}$.
\end{prp}

\begin{proof} Applying Theorem \ref{thDK} with a trivial
1-dimensional Virasoro module, we see that $\mathcal M$ as a
module for the Lie algebra $\Venone$.  By restriction, view
$\mathcal M$ as an $\tsl_{N+1}$-module. The top of the module
$\mathcal M$ is the tensor module
$$T (W, \gamma) = q^\gamma \C [q_1^{\pm 1}, \ldots, q_N^{\pm 1}]  \otimes W .$$
To show that $\mathcal M$ is a generalized Wakimoto module for
$\tsl_{N+1}$, we need to compare the characters of $\C[u_{pj},
v_{pj}] \otimes U(\hgl^-_{N})$ and $U(\tsl_{N+1}^-)$. We have
$$\chara \C[u_{pj}, v_{pj}] \times \chara U(\hgl^-_{N}) = \prod_{k=1}^\infty (1-s^k)^{-2N} \times
\prod_{k=1}^\infty (1-s^k)^{-N^2} = \chara U(\tsl_{N+1}^-) .$$
This completes the proof of the proposition.

\end{proof}

Theorem \ref{thDK} describes irreducible $\D \ltimes \K$-modules.
We would like to study their reductions to subalgebra $\D$. In
general, when reduced to a subalgebra, modules become reducible.
Here, however, the link with generalized Wakimoto modules for
$\tsl_{N+1}$ and the result of \cite{BMZ} for $N=1$, give us hope
that the situation may be better than one would expect a priory.

\

\section{Duality for modules over the Lie algebra of vector
fields}

In this section we will establish a duality for the class of modules described in Theorem \ref{thDK} (ii), that will be useful for
the study of their irreducibility as modules over Lie algebra $\D$. We begin by looking at this question in a general setup.

Let $\mathcal L$ be a $\Z$-graded Lie algebra $\mathcal L =
\mathop\oplus\limits_{n\in\Z} \mathcal L_n$ with an
anti-involution $\sigma$ such that $\sigma(\mathcal L_n) =
\mathcal L_{-n}$. Extend $\sigma$ to the universal enveloping
algebra $U(\mathcal L)$ by $\sigma(ab) = \sigma(b) \sigma(a)$. Let
$\mathcal L_{\pm} = \mathop\oplus\limits_{n>0} \mathcal L_{\pm
n}$.
 Suppose $T_1$,
$T_2$ be two $\mathcal L_0$-modules with a bilinear pairing
$$T_1  \times T_2 \rightarrow \C \eqno{(\bilpair)}$$
such that
$$ \left< x w_1 , w_2 \right> = \left< w_1, \sigma (x) w_2 \right> $$
for all $x \in \mathcal L_0$, $w_1 \in T_1$, $w_2 \in T_2$.

For an $\mathcal L_0$-module $T$ we let $\mathcal L_+$ act on $T$
trivially and construct the generalized Verma modules for
$\mathcal L$: $M(T) = U(\mathcal L_-) \otimes T$.

The generalized Verma module $M(T)$ inherits the $\Z$-grading (assuming degree of $T$ to be zero).
 Define the {\it radical} of a generalized Verma module $M(T)$ as the maximal homogeneous submodule trivially intersecting with
the top $T$. If $T$ is an irreducible $\mathcal L_0$-module then
the quotient $L(T)$ of $M(T)$ by its radical is an irreducible
module for $\mathcal L$.

Consider the Shapovalov projection
$$ S: \, U(\mathcal L) \rightarrow U(\mathcal L_0) $$
with kernel $\mathcal L_- U(\mathcal L) + U(\mathcal L) \mathcal
L_+$.

 Define a bilinear pairing
$$ M(T_1) \times M(T_2) \rightarrow \C,$$
defined by
$$ \langle a  w_1, b w_2\rangle = \langle w_1, S(\sigma(a) b) w_2\rangle \eqno{(\copair)}$$
for $a, b \in U(\mathcal L_-)$, $w_1 \in T_1$, $w_2 \in T_2$.

\

\begin{prp}\label{contr} (i) The pairing (\copair) is {\it contragredient},
i.e.,
$$ \langle x u, v \rangle = \langle  u, \sigma(x) v\rangle $$
for all $x\in \mathcal L$, $u\in M(T_1)$, $v\in M(T_2)$.

(ii) If $n \neq k$ then $ \langle M(T_1)_n, M(T_2)_k \rangle =0$.

(iii) The radicals of $M(T_1)$, $M(T_2)$ are in the kernel of the pairing.

(iv) Assume that $T_1$, $T_2$ are irreducible $\mathcal
L_0$-modules and the pairing (\bilpair) is non-zero. Then the
contragredient pairing factors to the simple modules
$$ L(T_1)  \times L(T_2)   \rightarrow \C ,$$
on which it is non-degenerate.
\end{prp}

This proposition is standard (cf., \cite{MP}, Proposition 2.8.1)
and its proof is left to the reader as an exercise.

Next we will apply this proposition to establish the duality for
the bounded modules described in Theorem \ref{thDK} (ii).

 We consider the following anti-involution on $\D \ltimes \K$:
$$\sigma (t_0^j t^r d_a ) = t_0^{-j} t^{-r} d_a, \,\quad
\sigma (t_0^j t^r k_a ) = t_0^{-j} t^{-r} k_a, \,\quad  a = 0, \ldots, N.$$

For a finite-dimensional simple $\gl_N$-module $W$, on which the identity matrix $I$
acts as multiplication by $\alpha \in\C$,
let $W^*$ be the dual space to $W$ with $\sln_N$-module structure of the dual module, but
with $I$ acting as scalar $N - \alpha$. The natural pairing between $W$ and $W^*$ satisfies
$$ \langle E^{ab} w | w^* \rangle = \langle w | -E^{ab} w^* + \delta_{ab} w^* \rangle .$$

\begin{thm} \label{duality} There exists a non-degenerate contragredient
pairing of simple $\D \ltimes \K$-modules defined in Theorem
\ref{thDK} (ii):
$$L(W,\gamma,h) \times L(W^*,\gamma,h) \rightarrow \C, \eqno{(\bigpair)}$$
satisfying
$$\langle x u, v \rangle = \langle u, \sigma(x) v \rangle,$$
for all $x \in \D \ltimes \K$, $u\in L(W,\gamma,h)$, $v\in  L(W^*,\gamma,h)$.
\end{thm}

For the proof of this theorem we will use an alternative construction
of the simple $\D \ltimes \K$-module $L(W,\gamma,h)$,
which is discussed in \cite{B}. These modules may be abstractly
defined using approach of Theorem \ref{TBB}. The {\it top} of the
module $L(W,\gamma,h)$ is the space
$$T = T(W,\gamma,h) = q^\gamma \C [q_1^{\pm 1}, \ldots, q_N^{\pm 1}] \otimes W \otimes \vh ,$$ which is a module for
the zero degree component $\D_0 \ltimes \K_0$ of $\D \ltimes \K$ with respect to its $\Z$-grading by degree in $t_0$.
The positive part of $\D \ltimes \K$ acts on $T(W,\gamma,h)$ trivially, and we can consider the induced
$\D \ltimes \K$-module $M(T)$. The induced module has a unique irreducible quotient, which is
isomorphic to $L(W,\gamma,h)$.

 The action of $\D_0 \ltimes \K_0$ on $T(W,\gamma,h)$ can be derived from Theorem \ref{thDK} (i) (see Theorem 6.4 in \cite{B} for details) and $T$ is essentially a tensor module that we discussed above:
$$ t^r k_0 (q^\mu \otimes w \otimes \vh) = q^{\mu+r} \otimes w \otimes \vh, \eqno{(\contra)}$$
$$ t^r k_a (q^\mu \otimes w \otimes \vh) = 0,\eqno{(\contrb)}$$
$$ t^r d_a (q^\mu \otimes w \otimes \vh) =
\mu_a q^{\mu+r} \otimes w \otimes \vh + \sum_{p=1}^N r_p q^{\mu+r} \otimes E^{pa} w \otimes \vh, \eqno{(\contrc)}$$
$$ t^r d_0 (q^\mu \otimes w \otimes \vh) = \beta \; q^{\mu+r} \otimes w \otimes \vh, \eqno{(\contrd)}$$
where
$$\beta = -h - {{\Omega_W} \over {2(N+1)}} - {{\alpha (\alpha - N)} \over {2N}}, \eqno{(\betavalue)}$$
$a = 1, \ldots, N$, $r\in \Z^N$, $\mu \in \gamma + \Z^N$. Here $\Omega_W$ is the scalar
with which the Casimir operator of $\sln_N$ acts on $W$.

 The claim of Theorem \ref{duality}  follows from Proposition
 \ref{contr} and the following lemma.

\begin{lem} \label{Tpair} Let $W$ be a simple finite-dimensional
$\gl_N$-module on which the identity matrix acts as scalar
$\alpha$, and let $W^*$ be a $\gl_N$-module which is dual to $W$
as a $\sln_N$-module, and on which the identity matrix acts as
scalar $N-\alpha$.
 Then the pairing
$$ T(W,\gamma,h) \times T(W^*,\gamma,h) \rightarrow \C,$$
given by
$$\langle q^\mu \otimes w \otimes \vh | q^\eta \otimes w^* \otimes \vh \rangle = \delta_{\mu,\eta} \langle w | w^* \rangle$$
is a non-degenerate contragredient pairing of $\D_0 \ltimes \K_0$-modules.
\end{lem}

\begin{proof} Using (\contra) we get
$$\langle t^r k_0 (q^\mu \otimes w \otimes \vh) | q^\eta \otimes w^* \otimes \vh \rangle
= \delta_{\mu+r,\eta} \langle w | w^* \rangle$$
$$ = \delta_{\mu,\eta-r} \langle w | w^* \rangle
= \langle q^\mu \otimes w \otimes \vh | \sigma(t^r k_0) (q^\eta \otimes w^* \otimes \vh) \rangle.$$
In case of $t^r k_a$, $a=1,\ldots, N$, both left and right hand sides are zero.
Let us verify the contragredient property for the action of  $t^r d_a$, $a=1,\ldots, N$:
$$\langle t^r d_a (q^\mu \otimes w \otimes \vh) | q^\eta \otimes w^* \otimes \vh \rangle $$
$$= \langle \mu_a q^{\mu+r} \otimes w \otimes \vh + \sum_{p=1}^N r_p q^{\mu+r} \otimes E^{pa} w \otimes \vh
| q^\eta \otimes w^* \otimes \vh \rangle $$
$$= \delta_{\mu+r,\eta} \left( \mu_a \langle w | w^* \rangle + \sum_{p=1}^N r_p \langle E^{pa} w | w^* \rangle \right)$$
$$= \delta_{\mu,\eta-r} \left(\mu_a \langle w | w^* \rangle + \sum_{p=1}^N r_p \langle w | -E^{pa} w^* + \delta_{pa} w^* \rangle \right)$$
$$= \delta_{\mu,\eta-r} \left( (\eta_a - r_a) \langle w | w^* \rangle - \sum_{p=1}^N r_p \langle w | E^{pa} w^* \rangle
 + r_a \langle w | w^* \rangle \right)$$
$$= \delta_{\mu,\eta-r} \left(\eta_a  \langle w | w^* \rangle - \sum_{p=1}^N r_p \langle w | E^{pa} w^* \rangle \right)$$
$$= \langle q^\mu \otimes w \otimes \vh | \sigma(t^r d_a) (q^\eta \otimes w^* \otimes \vh) \rangle .$$
Finally, to check the case of $t^r d_0$, we note that the constant
$\beta$ in (\betavalue) is the same for $T(W,\gamma,h)$ and
$T(W^*,\gamma,h)$. This follows from the fact that the Casimir
operator for $\sln_N$ acts with the same scalar on $W$ and $W^*$,
while the last term in (\betavalue) is invariant under the
substitution $\alpha \mapsto N - \alpha$. Thus the computation in
the case of $t^r d_0$ is analogous to the case of $t^r k_0$. This
completes the proof of the lemma.
\end{proof}
\

\begin{rem} The pairing (\bigpair) is in fact a product of
contragredient pairings of tensor factors
$$ \Mhyp(\gamma) \times \Mhyp(\gamma) \rightarrow \C, \quad
\Lgl (W) \times \Lgl (W^*) \rightarrow \C, \quad
\Lvir(h) \times \Lvir(h) \rightarrow \C, $$
with respect to appropriate anti-involutions of corresponding Lie algebras.

\end{rem}

\begin{rem} The duality of Theorem \ref{duality} can be alternatively
constructed via vertex algebra approach, using the definition of
the contragredient module over a vertex algebra (see section 5.2
in \cite{FHL}).
\end{rem}

 One of the goals of this paper is to analyze which of the modules defined in Theorem \ref{thDK} (ii)
remain irreducible after restriction to $\D$. First of all, let us look at the question of irreducibility
of the top $T(W,\gamma,h)$ as a module over $\D_0$.

\begin{lem} \label{irrtop} Let $W$ be an irreducible finite-dimensional
$\gl_N$-module. The module $T(W,\gamma,h)$ is reducible as a
$\D_0$-module if and only if it is reducible as a $\VenN$-module
(see Theorem \ref{DRc}) and $h=0$.
\end{lem}

\begin{proof} Clearly, if $T(W,\gamma,h)$ is reducible as a
$\D_0$-module, it must also be reducible as a $\VenN$-module. By
Theorem \ref{DRc}, all such modules appear in the de Rham complex.
Note that
$$\D_0 = \VenN \oplus \C[t_1^{\pm 1}, \ldots, t_N^{\pm 1}] d_0,$$
and by (\contrd), $t^r d_0$ acts on $T(W,\gamma,h)$ as
multiplication by $\beta q^r$. It is well known that in the
modules of differential forms there are no proper submodules that
are $\C[q_1^{\pm 1}, \ldots, q_N^{\pm 1}]$-invariant. Thus for
$T(W,\gamma,h)$ to be reducible as a $\D_0$-module it is necessary
and sufficient that it is reducible as a $\VenN$-module and the
value of $\beta$ given by (\betavalue) is zero. Let us analyze the
values of $\beta$ for the modules in the de Rham complex. For the
modules $\Omega^0 (\T^N)$ and $\Omega^N (\T^N)$ the
$\sln_N$-module $W$ is trivial, so the Casimir operator acts with
constant $\Omega_W = 0$, while the identity matrix acts on $W$
with scalars $\alpha=0$ and $\alpha=N$ respectively. Simplifying
the expression in (\betavalue) we get in this case that
$\beta=-h$. In case of the modules of $k$-forms $\Omega^k (\T^N)$,
$k=1,\ldots,N-1$, the highest weight of the corresponding
$\sln_N$-module $W$ is the fundamental weight $\omega_k$. A
standard computation shows that in this case the Casimir operator
acts with the scalar
$$\Omega_W = (\omega_k | \omega_k + 2\rho) = k (N-k) {{N+1} \over N} . \eqno{(\casim)}$$
Since the identity matrix acts with the scalar $\alpha = k$, the formula (\betavalue)
again simplifies to $\beta = -h$. This implies the claim of the Lemma.

\end{proof}

Consider now an irreducible $\D \ltimes \K$ module $L(W,\gamma,h)$ described
in Theorem \ref{thDK} (ii), and assume that its top $T(W,\gamma,h)$ is irreducible
as a $\D_0$-module. To show that $L(W,\gamma,h)$ remains irreducible as a
module over $\D$, it is sufficient to establish two properties:

\vskip 0.3cm
(C) Every critical vector of $L(W,\gamma,h)$ (i.e., annihilated by $\D_+$)
belongs to its top $T(W,\gamma,h)$.

\vskip 0.3cm
(G) $L(W,\gamma,h)$ is generated by its top $T(W,\gamma,h)$ as a module over $\D_-$.

\

 The following standard observation will be quite useful:

\begin{lem} \label{CGduality} Condition (C) holds for the module
$L(W,\gamma,h)$ if and only if condition (G) holds for
$L(W^*,\gamma,h)$.
\end{lem}

\ifnum \value{version}=\value{short} 
This lemma is an immediate consequence of Theorem \ref{duality}
and we omit its proof.
\else
 \begin{proof} We use the existence of a non-degenerate
contragredient pairing of $\D \ltimes \K$-modules:
$$ L(W,\gamma,h) \times L(W^*,\gamma,h) \rightarrow \C .$$
If $L(W,\gamma,h)$ has a vector annihilated by $\D_+$, which does not belong
to the top, it also has a homogeneous vector with this property. Suppose $u$ is
a critical vector of degree $s\in\Z^{N+1}$. Since the pairing is non-degenerate,
there exists a vector $v$ in $L(W^*,\gamma,h)$ of the same degree, such that
$\langle u | v \rangle \neq 0$. If property (G) holds for $L(W^*,\gamma,h)$,
$v$ can be written as $v = \sum_i x_i v_i$, where $x_i \in \D_-, v_i \in L(W^*,\gamma,h)$.
Applying the contragredient property, we get
$$ \langle u | v \rangle = \langle u | \sum_i x_i v_i \rangle =
\sum_i \langle \sigma(x_i) u | v_i \rangle .$$
The last expression is zero since $\sigma(x_i) \in \D_+$ and $u$ is a critical vector.
This gives a contradiction, which implies that property (G) does not hold for $L(W^*,\gamma,h)$.

To prove the converse, assume that the component of degree $s$ in $L(W^*,\gamma,h)$ is not
generated by $\D_-$ acting on the top. Let $V$ be the intersection of that homogeneous
component with the space $U(\D_-) T(W^*,\gamma,h)$. Since the pairing is non-degenerate, we
can find a non-zero vector $u$ in the degree $s$ component of $L(W,\gamma,h)$, such that
$ \langle u | V \rangle = 0$. If property (C) holds for the module $L(W,\gamma,h)$,
there exists a homogeneous $y \in U(\D_+)$, such that $z = y u$ is a non-zero vector in $T(W,\gamma,h)$.
Let $z^\prime \in T(W,\gamma,h)$ be such that $\langle z | z^\prime \rangle \neq 0$.
Then
$$ \langle z | z^\prime \rangle = \langle y u | z^\prime \rangle = \langle u | \sigma(y) z^\prime \rangle .$$
But $\sigma(y) \in U(\D_-)$, thus $\sigma(y) z^\prime \in V$,
which leads to a contradiction. The lemma is now proved.

\end{proof}
\fi

Lemma \ref{CGduality} reduces the question of irreducibility of the family of
modules $L(W,\gamma,h)$ to the question of existence of critical
vectors. If both $L(W,\gamma,h)$ and $L(W^*,\gamma,h)$ satisfy
condition (C), then they are both irreducible as modules over
$\D$.

\

\

\section{Critical vectors}

In this section we will establish a necessary condition for the
existence of non-trivial critical vectors in the modules
$L(W,\gamma,h)$, which together with Lemma \ref{CGduality} will
give a sufficient condition for the irreducibility of such
modules.

\begin{thm} \label{crith} Let $W$ be an irreducible finite-dimensional
$\gl_N$-module, $\gamma\in \C^N, h\in\C$. Every critical vector
(i.e., annihilated by $\D_+$) in the module $L(W,\gamma,h)$
belongs to its top, unless $h=0$ and $W$ is either a trivial
$\sln_N$-module with identity matrix acting with scalar $\alpha =
N - m N$, $m = 1,2, \ldots$, or the highest weight of $W$ is a
fundamental weight $\omega_k$, $k=1,\ldots,N-1$, with identity
matrix acting by scalar $\alpha = k - mN$, $m = 1,2, \ldots$.

\end{thm}

 We will call a module $L(W,\gamma,h)$ {\it exceptional} if $h=0$,
the identity matrix acts on $W$ by an integer $k \in\Z$ and $W$ is a trivial
one-dimensional $\sln_N$-module when $k = 0 \mod N$ or has a fundamental
highest weight $\omega_{k^\prime}$ with $1 \leq k^\prime \leq N-1$ and
$k = k^\prime \mod N$.

\begin{thm} \label{irr}  Let $W$ be an irreducible finite-dimensional
$\gl_N$-module, $\gamma\in \C^N, h\in\C$. Every non-exceptional module 
$L(W,\gamma,h)$ is irreducible as a $\VenN$-module.
\end{thm}

 Theorem \ref{irr} is an immediate consequence of Theorem \ref{crith} and Lemma 
\ref{CGduality}.  
The proof of Theorem \ref{crith} will be split into a sequence of lemmas.

\

\begin{lem} \label{cone} Let $g \in L(W,\gamma,h)$ be a critical
vector. Then $g$ does not depend on variables $\{ u_{pj} \}$,
i.e., $g$ belongs to the subspace
$$q^\gamma \C[q_1^{\pm 1}, \ldots, q_N^{\pm 1}] \otimes \C [v_{pj} |^{p=1,\ldots,N}_{j=1,2,\ldots}]
\otimes \Lgl (W) \otimes \Lvir (h) . \eqno{(\indu)}$$
\end{lem}

\begin{proof} The algebra $\D_+$ contains the elements $t_0^j d_p$ with
${p=1,\ldots,N}$, $j \geq 1$, which act as $\partial \over
{\partial u_{pj}}$. The condition $(t_0^j d_p) g = 0$ implies the
claim of the lemma.

\end{proof}

  For a formal series $a(z)$ we denote by $a(z)_-$ its part
that only involves negative powers of $z$.
 Recalling that $d_a(r,z) = \sum\limits_{j\in\Z} t_0^j t^r d_a z^{-j-1}$, we have $(z d_a(r,z))_-  g = 0$.
Using (\actthree) for the action of $d_a(r,z)$ and taking into account that $g$ does not depend on
$\{ u_{pj} \}$, we get
$$q^r \left( \exp \left( \sum_{p=1}^N r_p \sum_{j=1}^\infty u_{pj} z^j \right)
\left( \sum_{i=1}^\infty i v_{ai} z^i + q_a {\partial \over
\partial q_a} + \sum_{p=1}^N r_p \sum_{k \in \Z} E^{pa}_{(k)}
z^{-k} \right) \right.$$
$$ \left. \times \exp \left( - \sum_{p=1}^N r_p \sum_{j=1}^\infty {z^{-j} \over j}
{\partial \over \partial v_{pj}} \right) \right)_- g = 0 .$$

Let us project to the subspace (\indu), setting $u_{pj} = 0$ in the above equality. Also, since the operator
of multiplication by $q^r$ is invertible, we can drop it. We then get
$$ P_a(r,z)_- g = 0, \quad a = 1, \ldots, N,  \eqno{(\simpg)}$$
where
$$P_a(r,z) =
\left( \sum_{i=1}^\infty i v_{ai} z^i + q_a {\partial \over
\partial q_a} + \sum_{p=1}^N r_p \sum_{k \in \Z} E^{pa}_{(k)}
z^{-k} \right) \exp \left( - \sum_{p=1}^N r_p \sum_{j=1}^\infty
{z^{-j} \over j} {\partial \over \partial v_{pj}} \right). $$

 At this point we find it convenient to make a change of variables $x_{aj} = j v_{aj}$. In these notations $P_a(r,z)$
takes form
$$P_a(r,z) =
\left( \sum_{i=1}^\infty x_{ai} z^i+ q_a {\partial \over \partial
q_a} + \sum_{p=1}^N r_p \sum_{k \in \Z} E^{pa}_{(k)} z^{-k}
\right) \exp \left( - \sum_{p=1}^N r_p \sum_{j=1}^\infty {z^{-j} }
{\partial \over \partial x_{pj}} \right). \eqno{(\Prz)} $$

Let us expand $P_a(r,z)$ in a formal series in variables $r = (r_1, \ldots, r_N)$:
$$ P_a(r,z) = \sum_{s \in \Z_+^N} r^s P_{as} (z) .$$
It is easy to see that for any $j\in\Z$ and any vector $g^\prime$, there are only finitely many $s\in \Z_+^N$ such
that the coefficient at $z^j$ in $P_{as} (z) g^\prime$ is non-zero.
Thus the coefficient at $z^j$ in $\sum\limits_{s \in \Z_+^N} r^s P_{as} (z) g$ is a polynomial in $r$.
Since for each $j<0$ these polynomials vanish for all $r \in \Z^N$, we conclude that for all $s \in \Z_+^N$, $a=1, \ldots, N$,
$$P_{as}(z)_- g = 0.$$
Note that for $s=0$ this equation is trivial. Let us consider the case $s \in \Z_+^N$, with $s_p =1$ and $s_i = 0$ for $i\neq p$. This gives us an equation
$$ \sum_{k=1}^\infty z^{-k} E^{pa}_{(k)}  g = \left( \left( \sum_{i=1}^\infty x_{ai} z^i +
q_a {\partial \over \partial q_a} \right) \left( \sum_{k=1}^\infty
{z^{-k} } {\partial \over \partial x_{pk}} \right) \right)_- g $$
$$ = \left( \sum_{k=1}^\infty \left( \sum_{i=1}^{k-1} x_{ai} z^i +
q_a {\partial \over \partial q_a} \right) {z^{-k} } {\partial
\over \partial x_{pk}} \right) g. \eqno{(\Eqq)} $$

 Substituting (\Eqq) into (\Prz) we get
$$P_a^\prime (r,z)_- g = 0,   \eqno{(\Pz)} $$
where
$$P_a^\prime (r,z) = \left( \sum_{i=1}^\infty x_{ai} z^i +
q_a {\partial \over \partial q_a} + \sum_{p=1}^N r_p
\sum_{k=0}^\infty E^{pa}_{(-k)} z^k  \right) \exp \left( -
\sum_{p=1}^N r_p \sum_{j=1}^\infty {z^{-j} } {\partial \over
\partial x_{pj}} \right)  $$
$$+ \exp \left( - \sum_{p=1}^N r_p \sum_{j=1}^\infty {z^{-j} }
{\partial \over \partial x_{pj}} \right) \left( \sum_{p=1}^N r_p
\sum_{k=1}^\infty \left( \sum_{i=1}^{k-1} x_{ai} z^i + q_a
{\partial \over \partial q_a} \right) {z^{-k} } {\partial \over
\partial x_{pk}} \right) .
$$

 Our module is $\Z^{N+1}$-graded via the action of operators $d_0, \ldots, d_N$, and without loss of generality we may assume that $g$ is homogeneous relative to $\Z^{N+1}$-grading. We will call the eigenvalue of $\beta \, \Id -d_0$
the {\it degree} of $g$.  We use the negative sign here to make the degree non-negative. In fact the $\Z$-grading by
degree may be defined on each of the tensor factors $\C[x_{pj} |^{p=1, \ldots, N}_{j=1,2,\ldots}]$, $\Lgl$ and $\Lvir$
by
$$\deg(x_{pj}) = \deg(E^{ab}_{(-j)}) = \deg(L_{(-j)}) = j,$$
$$\deg(W) = \deg(\vh) = 0.$$

On the space $\C[x_{pj} |^{p=1, \ldots, N}_{j=1,2,\ldots}]$ we will also consider a refinment of $\Z$-grading by degree,
where we will compute the degree in each of the $N$ families of variables. For each $a = 1, \ldots, N$, define
$$\deg_a (x_{pj}) = j \delta_{ap} .$$
Then for a monomial $y \in \C[x_{pj} |^{p=1, \ldots, N}_{j=1,2,\ldots}]$ we have
$$\deg (y) = \sum_{a=1}^N \deg_a (y) .$$

In addition to the degree of monomials, we will consider another $\Z^N$-grading by {\it length}, where
$$\len_a (x_{pj}) =  \delta_{ap} $$
and define the {\it total length} to be
$$\len (y) = \sum_{a=1}^N \len_a (y) .$$

 Let us fix homogeneous bases $\{ y_i^\prime \}$, $\{ y_j^\pprime \}$, $\{ y_k^\ppprime \}$ in the spaces
$\C[x_{pj} |^{p=1, \ldots, N}_{j=1,2,\ldots}]$, $\Lgl$ and $\Lvir$ respectively.
Then we can expand $g$ into a finite sum
$$ g = \sum_{ijk} \alpha_{ijk} q^\mu y_i^\prime \otimes y_j^\pprime \otimes y_k^\ppprime . \eqno{(\gijk)}$$
Note that in the above decomposition $\deg(g) = \deg(y_i^\prime)+\deg(y_j^\pprime)+\deg(y_k^\ppprime)$.
Since equations (\simpg), (\Eqq) and (\Pz) do not involve any operators acting on the component $\Lvir$,
we conclude that these must be satisfied not only by $g$, but also by each of the components
$$ g_k = \sum_{ij} \alpha_{ijk} q^\mu y_i^\prime \otimes y_j^\pprime \otimes y_k^\ppprime .$$

\begin{lem} \label{topvect}  Let $g$ be a homogeneous non-zero
critical vector. Then in the decomposition (\gijk) there exist
$y_j^\pprime \in W$ and $y_k^\ppprime \in \C \vh$ with
$\alpha_{ijk} \neq 0$ for some $i$.
\end{lem}

\begin{proof}  Let us rewrite (\gijk) as
$$g =  \sum_{ij}  q^\mu y_i^\prime \otimes y_j^\pprime \otimes \left( \sum_k \alpha_{ijk} y_k^\ppprime \right) . $$
Consider the smallest degree $n_0$ of $y_k^\ppprime$ for which $\alpha_{ijk} \neq 0$ for some $i, j$. We claim that $n_0 = 0$.
Otherwise, since $\Lvir$ is irreducible, there exists a raising operator $\omega^\vir_{(n)}$, $n \geq 2$, such that
$$ \omega^\vir_{(n)}  \sum_{k \atop \sdeg (y_k^\ppprime) = n_0} \alpha_{ijk} y_k^\ppprime \neq 0 .$$
Consider now the Virasoro operator acting on all
three factors of the tensor product: $\omega_{(n)} = \omega^\hyp_{(n)} + \omega^{\gl_N}_{(n)} + \omega^\vir_{(n)}$. The part
of $\omega_{(n)} g$ involving the terms of the smallest degree in the component $\Lvir$ will be
$$ \sum_{ij}  q^\mu y_i^\prime \otimes y_j^\pprime \otimes
\omega^\vir_{(n)} \sum_{k \atop \sdeg (y_k^\ppprime) = n_0}  \alpha_{ijk} y_k^\ppprime \neq 0 . $$
Thus $\omega_{(n)} g \neq 0$. However  operator $\omega_{(n)}$ represents $-t_0^{n-1}d_0 \in \D_+$ and must annihilate
$g$ since $g$ is a critical vector. This is a contradiction, which implies that $n_0 = 0$.

Let  $\tg =  q^\mu \sum_{ij} \alpha_{ij}  y_i^\prime \otimes y_j^\pprime \otimes \vh $ be the projection of $g$ to the space
$$q^\mu \C [x_{pj} |^{p=1,\ldots,N}_{j=1,2,\ldots}] \otimes \Lgl (W) \otimes \vh .$$
We just proved that $\tg \neq 0$, and it was noted earlier that it satisfies equation (\Eqq).
Let $n_1$ be the smallest degree of  $y_j^\pprime$ such that $\alpha_{ij} \neq 0$ for some $i$. To complete the proof of the lemma,
we need to show that $n_1 = 0$. If $n_1 > 0$ using the same argument as above we see that there exists a raising operator $E^{pa}_{(n)}$,
$n \geq 1$ such that $E^{pa}_{(n)} \tg$ will have a non-zero component with terms in $\Lgl (W)$ of degree $n_1 - n$. However
 the equation (\Eqq) implies that all factors from $\Lgl (W)$ that appear in $E^{pa}_{(n)} \tg$ have degrees at least $n_1$. This contradiction
implies $n_1 = 0$, and the lemma is proved.

\end{proof}

Let $\gb$ be the projection of the critical vector $g$ to the space
$$q^\mu \C [x_{pj} |^{p=1,\ldots,N}_{j=1,2,\ldots}] \otimes W \otimes \vh . \eqno{(\VW)}$$
By the above Lemma, $\gb \neq 0$.  Let us take the projection of the equation (\Pz)
to the space (\VW). This yields
$$P_a^\pprime (r,z)_- \gb = 0, \eqno{(\Pza)}$$
where
$$P_a^\pprime (r,z) = \left( \sum_{i=1}^\infty x_{ai} z^i +
q_a {\partial \over \partial q_a} + \sum_{p=1}^N r_p E^{pa}_{(0)}
\right) \exp \left( - \sum_{p=1}^N r_p \sum_{j=1}^\infty {z^{-j} }
{\partial \over \partial x_{pj}} \right)  $$
$$+ \exp \left( - \sum_{p=1}^N r_p \sum_{j=1}^\infty {z^{-j} }
{\partial \over \partial x_{pj}} \right) \left( \sum_{p=1}^N r_p
\sum_{k=1}^\infty \left( \sum_{i=1}^{k-1} x_{ai} z^i + q_a
{\partial \over \partial q_a} \right) {z^{-k} } {\partial \over
\partial x_{pk}} \right) .
$$

Again we decompose $P_a^\pprime (r,z)$ in a formal power series in $r$, $P_a^\pprime (r,z) = \sum\limits_{s\in\Z^N_+} r^s P^\pprime_{as} (z)$,
and we have $P^\pprime_{as} (z)_- \gb = 0$ for all $s \in \Z^N_+$.

Next we will consider the grading of $\gb$ by (total) length. Note
that the operator $P^\pprime_{as} (z)$ has two homogeneous
components with respect to the length grading -- one that
decreases the length by $s_1 + \ldots + s_N - 1$ and the other
(that contains terms involving $q_a {\partial \over \partial
q_a}$) by $s_1 + \ldots + s_N$. Let $\gh$ be the maximum length
component of $\gb$ and suppose $\len(\gh) = \ell$. Denote by
$Q_{as} (z)$ the component of $P^\pprime_{as} (z)$ that reduces
the length by $s_1 + \ldots + s_N - 1$. Then the component of
length $\ell + 1 - s_1 - \ldots -s_N$ in $P^\pprime_{as} (z)_-
\gb$ is $Q_{as} (z)_- \gh$. Thus
$$Q_{as} (z)_- \gh = 0. \eqno{(\Pzb)} $$
Assembling back the generating series $Q_a (r,z) =  \sum\limits_{s\in\Z^N_+} r^s Q_{as} (z)$, we get
$$Q_a (r,z) = \left( \sum_{i=1}^\infty x_{ai} z^i  + \sum_{p=1}^N r_p  E^{pa}_{(0)} \right)
\exp \left( - \sum_{p=1}^N r_p \sum_{j=1}^\infty {z^{-j} }
{\partial \over \partial x_{pj}} \right)  $$
$$+ \exp \left( - \sum_{p=1}^N r_p \sum_{j=1}^\infty {z^{-j} } {\partial \over \partial x_{pj}} \right)
\left( \sum_{p=1}^N r_p \sum_{k=1}^\infty \left( \sum_{i=1}^{k-1}
x_{ai} z^i \right) {z^{-k} } {\partial \over \partial x_{pk}}
\right)
$$
$$= \left( \sum_{i=1}^\infty x_{ai} z^i  + \sum_{p=1}^N r_p  E^{pa}_{(0)} \right)
\exp \left( - \sum_{p=1}^N r_p \sum_{j=1}^\infty {z^{-j} }
{\partial \over \partial x_{pj}} \right)  $$
$$+ \left( \sum_{p=1}^N r_p \sum_{k=1}^\infty \left( \sum_{i=1}^{k-1} (x_{ai} - r_a z^{-i}) z^i \right)
{z^{-k} } {\partial \over \partial x_{pk}} \right)
 \exp \left( - \sum_{p=1}^N r_p \sum_{j=1}^\infty {z^{-j} } {\partial \over \partial x_{pj}} \right)$$
$$= \left( \sum_{i=1}^\infty x_{ai} z^i  \right)
\exp \left( - \sum_{p=1}^N r_p \sum_{j=1}^\infty {z^{-j} }
{\partial \over \partial x_{pj}} \right)  $$
$$+ \left( \sum_{p=1}^N r_p \sum_{k=1}^\infty \left( \sum_{i=1}^{k-1} x_{ai} z^i \right)
{z^{-k} } {\partial \over \partial x_{pk}} \right)
 \exp \left( - \sum_{p=1}^N r_p \sum_{j=1}^\infty {z^{-j} } {\partial \over \partial x_{pj}} \right)$$
$$+ \left(  \sum_{p=1}^N r_p  E^{pa}_{(0)} -  r_a \sum_{k=1}^\infty \sum_{p=1}^N (k-1)  r_p
{z^{-k} } {\partial \over \partial x_{pk}} \right)
 \exp \left( - \sum_{p=1}^N r_p \sum_{j=1}^\infty {z^{-j} } {\partial \over \partial x_{pj}} \right) . \eqno{(\oper)}$$

Our goal is to solve the system of equations (\Pzb) in the space  $\C [x_{pj} |^{p=1,\ldots,N}_{j=1,2,\ldots}] \otimes W$.
The solutions we are looking for are homogeneous in both degree and length. Let $\deg(\gh) = m$,
$\len(\gh) = \ell$. The equation (\Pzb) has trivial solutions with $\ell = 0$. These correspond to the critical vectors
in the top of the module. We are going to show that non-trivial solutions of (\Pzb) must have length $\ell =1$.

To establish this claim we will first analyze the equation (\Pzb) in cases $N=1$ and $N=2$. The case of the general $N$ will
follow from the following simple observation.
Consider a proper subset $S \subset \{1,\ldots,N \}$. Let us take a solution $\gh$ of (\Pzb) and specialize all variables $x_{pj}$ with
$p\not\in S$ to scalars, we will get a solution for (\Pzb) with a smaller $N$. To see this, set in (\Pzb) $r_p = 0$ for all $p \not\in S$
and restrict $a$  to the set $S$.
The information about the solutions of (\Pzb) with $N=1$ and $N=2$ may be used to establish properties
of solutions of this equation for a general $N$.

\

\begin{lem} \label{exq} Let $N=1$, and let $W$ be a 1-dimensional
$\gl_1$-module with identity matrix $I=E^{11}$ acting by scalar
$\alpha \in \C$. Let $\gh$ be a non-constant homogeneous (in both
length and degree) solution of (\Pzb) with $N=1$. Then $\len(\gh)
= 1$ and $\alpha = 1 -\deg(\gh)$.
\end{lem}

\begin{proof} First of all, let us rewrite (\oper) for the
case $N=1$. To simplify notations we will drop redundant indices
and write $x_i$ instead of $x_{pi}$, etc.,
$$Q(r,z) = \left( \sum_{i=1}^\infty x_{i} z^i
+  r \sum_{k=1}^\infty \left( \sum_{i=1}^{k-1} x_{i} z^i \right)
{z^{-k} } {\partial \over \partial x_{k}} + r \alpha -  r^2
\sum_{k=1}^\infty (k-1)  {z^{-k} } {\partial \over \partial x_{k}}
\right) \exp \left( -  r \sum_{j=1}^\infty {z^{-j} } {\partial
\over \partial x_{j}} \right) . \eqno{(\oner)}$$


 Let $\gh$ satisfy $Q(r,z)_- \gh = 0$. We may view $\gh$ as a polynomial of length $\ell > 0$
and degree $m \geq \ell$ in $\C [x_1, x_2, \ldots]$.
Choose a natural number $s$ such that
$$\left( \sum_{j=1}^\infty {z^{-j} } {\partial \over \partial x_{j}} \right)^{s+1} \gh = 0  \eqno{(\Rza)}$$
but
 $$ R(z) = \left( \sum_{j=1}^\infty {z^{-j} } {\partial \over \partial x_{j}} \right)^{s} \gh \neq 0.  \eqno{(\Rzb)}$$
Clearly $1 \leq s \leq \ell$.

Let us consider the coefficient at $r^{s+1}$ in the equation  $Q(r,z)_- \gh = 0$:
$$ {(-1)^s \over s!} \left(  \sum_{k=1}^\infty \left( \sum_{i=1}^{k-1} x_{i} z^i \right)
{z^{-k} } {\partial \over \partial x_{k}} \right)   R(z) + {(-1)^s
\over s!} \alpha R(z)$$
$$- {(-1)^{s-1} \over (s-1)!}
 \left(  \sum_{k=1}^\infty (k-1)  {z^{-k} } {\partial \over \partial x_{k}} \right)
 \left( \sum_{j=1}^\infty {z^{-j} } {\partial \over \partial x_{j}} \right)^{s-1} \gh = 0 .$$
Thus
$$ \left(  \sum_{k=1}^\infty \left( \sum_{i=1}^{k-1} x_{i} z^i \right)
{z^{-k} } {\partial \over \partial x_{k}} \right)   R(z)
 +  \alpha R(z) - z {d \over dz} R(z)  - s R(z) = 0. \eqno{(\Rzc)}$$
Consider an expansion $R(z) = R_n z^{-n} + R_{n+1} z^{-n-1} + \ldots + R_m z^{-m}$ with
$R_n \neq 0$. Let us look at the coefficient at $z^{-n}$ in the above equation:
$$\alpha R_n + n R_n - s R_n = 0,$$
which implies $\alpha = s - n \in \Z_-$.

Applying the operator $z \frac{d}{dz}$ to (\Rza) we get
$$\left( \sum_{k=1}^\infty k z^{-k} {\partial \over \partial x_k} \right) R(z) = 0. $$
Let us prove by induction that for all $j \geq 0$
$$\left( \sum_{k=1}^\infty k^j z^{-k} {\partial \over \partial x_k} \right) R(z) = 0. \eqno{(\Rzd)}$$
Suppose (\Rzd) holds for all $j^\prime \leq j$, $j \geq 1$.
Applying the operator $\sum_{k=1}^\infty k^j z^{-k} {\partial \over \partial x_k}$
to (\Rzc) and using the induction assumption, we get
$$ 0 = \left( \sum_{p=1}^\infty p^j z^{-p}  {\partial \over \partial x_p} \right)
\left(  \sum_{k=1}^\infty \left( \sum_{i=1}^{k-1} x_{i} z^i \right)
{z^{-k} } {\partial \over \partial x_{k}} \right)   R(z)
- \left( \sum_{k=1}^\infty k^j z^{-k} {\partial \over \partial x_k} \right)
z {d \over dz} R(z)$$
$$= \left(  \sum_{k=1}^\infty \left( \sum_{p=1}^{k-1} p^j  \right)
{z^{-k} } {\partial \over \partial x_{k}} \right)   R(z)
- \left( \sum_{k=1}^\infty k^{j+1} z^{-k} {\partial \over \partial x_k} \right) R(z)$$
$$= \left( {1 \over {j+1}} - 1 \right) 
\left( \sum_{k=1}^\infty k^{j+1} z^{-k} {\partial \over \partial x_k} \right) R(z),$$
which establishes the inductive step. In the above calculation we used the fact that
$\sum\limits_{p=1}^{k-1} p^j$ is a polynomial in $k$ with the leading term 
${k^{j+1} \over {j+1}}$.

Using the Vandermonde determinant argument we conclude from (\Rzd) that
$${\partial \over \partial x_{k}}  R(z)  = 0$$
for all $k$. This implies that all $R_i$ are scalars, but since they can't have equal degrees, we
conclude that $R_n$ is a non-zero scalar, while all other coefficients are zero. Without the loss of
generality we may thus assume
$$ R(z) = \left( \sum_{j=1}^\infty {z^{-j} } {\partial \over \partial x_{j}} \right)^{s} \gh = z^{-n}.$$
Thus $s = \len(\gh) = \ell$ and $n=\deg(\gh) = m$, and then $\alpha = \len(\gh) - \deg(\gh)$.
It remains to prove that $\len(\gh)=1$. We will reason by contradiction. Let us suppose that
$\ell = \len(\gh) > 1$ and consider
$$S(z) = \left( \sum_{j=1}^\infty {z^{-j} } {\partial \over \partial x_{j}} \right)^{\ell-1} \gh.$$
Since every term in $S(z)$ has length $1$, we can write
$$S(z) = \beta_1 x_1 z^{-m+1} + \beta_2 x_2 z^{-m+2} + \ldots + \beta_{m-\ell+1} x_{m-\ell+1} z^{-\ell+1} .$$
Now we look at the coefficient at $r^{\ell}$ in the equation  $Q(r,z)_- \gh = 0$:
$$0 =   {(-1)^\ell \over \ell!} \left( \left( \sum_{i=1}^\infty x_{i} z^i \right)
 \left( \sum_{j=1}^\infty {z^{-j} } {\partial \over \partial x_{j}} \right)^{\ell} \gh \right)_-$$
$$+  {(-1)^{\ell-1} \over (\ell-1)!} \left(  \sum_{k=1}^\infty \left( \sum_{i=1}^{k-1} x_{i} z^i \right)
{z^{-k} } {\partial \over \partial x_{k}} \right)
 \left( \sum_{j=1}^\infty {z^{-j} } {\partial \over \partial x_{j}} \right)^{\ell-1} \gh$$
$$-  {(-1)^{\ell-2} \over (\ell-2)!} \left(  \sum_{k=1}^\infty (k-1)  {z^{-k} } {\partial \over \partial x_{k}} \right)
 \left( \sum_{j=1}^\infty {z^{-j} } {\partial \over \partial x_{j}} \right)^{\ell-2} \gh$$
$$ +  {(-1)^{\ell-1} \over (\ell-1)!} \alpha   \left( \sum_{j=1}^\infty {z^{-j} } {\partial \over \partial x_{j}} \right)^{\ell-1} \gh .$$
Taking out the factor of $ {(-1)^\ell \over \ell!}$ we get
$$ 0 = \sum_{i=1}^{m-1} x_{i} z^{-m+i}
- \ell \left(  \sum_{k=1}^\infty \left( \sum_{i=1}^{k-1} x_{i} z^i
\right) {z^{-k} } {\partial \over \partial x_{k}} \right) S(z)$$
$$+ \ell z {d \over dz} S(z) + \ell(\ell-1) S(z) + \ell (m-\ell) S(z)$$
$$= \sum_{i=1}^{m-1} x_{i} z^{-m+i}
- \ell \sum_{k=1}^{m-\ell+1} \left( \sum_{i=1}^{k-1} x_i z^i \right) \beta_k z^{-m}$$
$$+\ell \sum_{k=1}^{m-\ell+1} (-m+k) \beta_k x_k z^{-m+k} + \ell(m-1) \sum_{k=1}^{m-\ell+1} \beta_k x_k z^{-m+k}$$
$$= \sum_{i=1}^{m-1} x_{i} z^{-m+i}
- \ell \sum_{i=1}^{m-\ell} \left( \sum_{k=i+1}^{m-\ell+1} \beta_k \right) x_i z^{-m+i}
+\ell \sum_{k=1}^{m-\ell+1} (k-1) \beta_k x_k z^{-m+k} .$$
We stress that the above calculation is only valid when $\ell > 1$. If $\ell > 2$ we immediately get a contradiction
since the coefficient at $z^{-1}$ yields $x_{m-1} = 0$. It only remains to rule out the case $\ell=2$.
In the latter case the equation simplifies to the following:
$$\sum_{k=1}^{m-1} x_{k} z^{-m+k}
- 2 \sum_{k=1}^{m-1} \left( \sum_{i=1}^{k-1} x_i z^i \right) \beta_k z^{-m}
+ 2 \sum_{k=1}^{m-1} (k-1) \beta_k x_k z^{-m+k} = 0. $$
Let us specialize this equation to $x_1 = \ldots = x_{m-1} = 1$, $z=1$. Then we get
$$(m-1) - 2 \sum_{k=1}^{m-1} (k-1) \beta_k + 2 \sum_{k=1}^{m-1} (k-1) \beta_k  = 0, $$
which implies $m=1$, which is impossible since $m = \deg(\gh)$ can not be less than $\ell = \len(\gh)$.
Thus $\ell = \len(\gh) = 1$ and the lemma is proved.

\end{proof}

 Next we are going to show that for a general $N$, a homogeneous non-trivial solution of (\Pzb) must have total
length $1$. The previous lemma implies that such a solution may have only two components with respect to
$\len_a$ grading for each $a$, where $\len_a$ may be either $0$ or $1$. To prove the general case, it is sufficient to consider $N=2$, since if a monomial has $\len_a + \len_b$ at most $1$ for any pair of distinct indices, then its
total length does not exceed $1$.

\

\begin{lem} \label{ctwo} Let $N=2$ and let $W$ be a finite-dimensional
$\gl_2$-module. Then any homogeneous (in both length and degree)
non-constant solution $\gh$ of (\Pzb) has total length $1$.
\end{lem}

\begin{proof} It is sufficient to consider the case of $W$
being irreducible since the equation (\Pzb) is compatible with the
$\gl_N$-module homomorphisms. Let us fix a basis \break $\{ w_n,
w_{n-1}, \ldots, w_{-n} \}$ of $W$, where $n \in {1\over 2} \Z_+$
and $(E^{11} - E^{22}) w_i = 2 i w_i$. Assuming that the identity
matrix acts on $W$ by scalar $\alpha$, we get
$$E^{11} w_i = \left({\alpha \over 2} + i\right) w_i \hbox{{\rm \ and \ }}
E^{22} w_i = \left({\alpha \over 2} -  i\right) w_i .$$

It follows from the previous lemma that every monomial in the
decomposition of $\gh$ has length at most $1$ with respect to each
of the two indices. Thus we only need to prove that $\gh$ can not
have total length $2$. We will reason by contradiction. If
$\len(\gh) = 2$ then for each monomial in $\gh$ both $\len_1$ and
$\len_2$ are $1$. Suppose $\deg(\gh) = m$. Let us write
$$\gh = \sum_{i=-n}^n \gh_i \otimes w_i .$$
By Lemma \ref{exq} we have
$$\deg_1 (\gh_i) = 1 -  \left({\alpha \over 2} + i\right), \quad
\deg_2 (\gh_i) = 1 -  \left({\alpha \over 2} - i\right),$$
so $m = \deg(\gh) = 2 - \alpha$. Let $b_i = \deg_1 (\gh_i) = {m \over 2} - i$ and
$c_i = \deg_2 (\gh_i) = {m \over 2} + i$. Then $\gh$ may be written as
$$\gh = \sum_{i=-n}^n \beta_i x_{1, b_i} x_{2,c_i} \otimes w_i , \eqno{(\ghtwo)}$$
where we set $\beta_i = 0$ whenever $b_i \leq 0$ or $c_i \leq 0$.

Let us take the equation derived from (\Pzb) by taking the coefficient at $r_1 r_2$ with $a=1$ and substitute
(\ghtwo) in it. We get
$$0 =  \sum_{i=-n}^n  \beta_i \left( \sum_{j=1}^{m-1}  z^{-m+j} x_{1,j}  \right) \otimes w_i
-  \sum_{i=-n}^n  \beta_i \left( \sum_{j=1}^{b_i-1}  z^{-m+j}  x_{1,j} \right)  \otimes w_i$$
$$-  \sum_{i=-n}^n  \beta_i \left( \sum_{j=1}^{c_i-1} z^{-m+j}  x_{1,j} \right)  \otimes w_i
-  \sum_{i=-n}^n \beta_i (c_i -1) z^{-m+b_i}  x_{1,b_i}  \otimes w_i$$
$$-  \sum_{i=-n}^n \beta_i (1-b_i) z^{-m+b_i}  x_{1,b_i} \otimes w_i
-  \sum_{i=-n}^n \beta_i  z^{-m+c_i}  x_{2,c_i} \otimes E^{21} w_i . \eqno{(\sumtwo)}$$

Note that only the last sum contains variables $x_{2,j}$. By equating this sum to zero, we conclude that
$\beta_i  E^{21} w_i = 0$ for all $i = -n, \ldots, n$. It follows that $\beta_i = 0$ for all $i \neq -n$.
Similarly, taking the same equation with $a=2$, we will get that $\beta_i = 0$ for all $i \neq n$. Thus the
only possibility for a non-zero solution is when $n=0$, which means that $W$ is 1-dimensional and $m$ is even.
In this case $b_0 = c_0 = {m \over 2}$ and
$\gh = x_{1, {m\over 2}} x_{2, {m\over 2}} \otimes w$, and the equation (\sumtwo) becomes
$$\sum_{j=1}^{m-1} x_{1,j} z^{-m+j}\otimes w - 2\sum_{j=1}^{ {m\over 2}-1} z^{-m+j} x_{1,j}  \otimes w
= 0,$$ which gives a contradiction. Thus the total length of $\gh$
must be $1$ and the lemma is proved.

\end{proof}

 Now we return to the general case. We proved that a non-trivial homogeneous solution $\gh$ of (\Prz) must have
length $1$. Suppose $\deg(\gh) = m$. Then $\gh$ can be written as
$$\gh = \sum_{p=1}^N x_{pm} \otimes w_p, \quad w_p \in W.$$
The equation (\Pzb) then simplifies as follows:
$$0 = \sum_{p=1}^N r_p \sum_{b=1}^N r_b E^{pa} w_b + r_a (m-1) \sum_{b=1}^N r_b w_b$$
$$ = \sum_{p=1}^N \sum_{b=1}^N r_p r_b \left( E^{pa} + \delta_{pa} (m-1) \right) w_b . \eqno{(\raq)}$$
Consider a new action $\rho^\prime$ of $\gln$ on $W$:
$$\rho^\prime (E^{pa})w = E^{pa}w + (m-1) \delta_{ap} w, \quad w \in W .$$
This gives the same structure of $W$ as an $\sln_N$-module, but now the identity matrix acts with
scalar $\alpha^\prime = \alpha + (m-1)N$. Then (\raq) is equivalent to the system of equations
$$\rho^\prime (E^{ca}) w_b + \rho^\prime (E^{ba}) w_c = 0 \eqno{(\skews)}$$
where $a, b, c = 1, \ldots, N$.

 We will also use a third $\gl_N$-action $\rho^\pprime$ on $W$:
$$\rho^\pprime (E^{pa})w = \rho^\prime (E^{pa})w + \delta_{ap} w = E^{pa}w + m \delta_{ap} w, \quad w \in W .$$
The identity matrix here acts with scalar $\alpha^\pprime = \alpha + mN$. For this action the equation (\skews)
may be written as
$$\rho^\pprime (E^{ca}) w_b + \rho^\pprime (E^{ba}) w_c = \delta_{ca} w_b + \delta_{ba} w_c . \eqno{(\skewtwo)}$$
We also have
$$  \sum_{p=1}^N \sum_{b=1}^N r_p r_b \rho^\pprime (E^{pa}) w_b
= r_a  \sum_{b=1}^N  r_b  w_b . \eqno{(\skewtri)}$$

 We are going to classify $\gln$-modules $W$ for which the system (\skews) has non-trivial solutions.
We will do this indirectly, linking this system with reducibility of tensor modules.

\begin{lem} \label{redten} Let $(W, \rho^\pprime)$ be a finite-dimensional
irreducible $\gl_N$-module. Let $\PP$ be the set of all solutions
$(w_1, \ldots, w_N) \in W \times \ldots \times W$ of the system of
equations (\skewtwo). Then the subspace
$$\tP = \left\{ \mathop\oplus\limits_{r \in \Z^N} q^r \otimes (r_1 w_1 + \ldots r_N w_N)
\big| (w_1, \ldots, w_N) \in \PP \right\}$$ is a $\VenN$-submodule
in the tensor module $T(W) = \C [q_1^{\pm 1}, \ldots, q_N^{\pm 1}]
\otimes W$, associated with $(W, \rho^\pprime)$.
\end{lem}

\begin{proof} Let $(w_1, \ldots, w_N) \in \PP$. Then using the
tensor module action and (\skewtri) we get
$$t^s d_a \left( q^r \otimes \sum_{b=1}^N r_b w_b \right)
= r_a  q^{r+s} \otimes \sum_{b=1}^N r_b w_b +
\sum_{p=1}^N \sum_{b=1}^N s_p r_b  q^{r+s} \otimes \rho^\pprime(E^{pa}) w_b$$
$$= \sum_{b=1}^N r_b q^{r+s} \otimes \sum_{p=1}^N (r_p + s_p)  \rho^\pprime(E^{pa}) w_b .$$
Fix $1 \leq a, i \leq N$. Set $\tw_p = \rho^\pprime(E^{pa}) w_i$,
$p=1,\ldots,N$. To complete the proof of the lemma, it is
sufficient to show that $(\tw_1, \ldots, \tw_N) \in \PP$. Instead
of working with (\skewtwo), it will be easier to check an
equivalent condition (\skews). Note that  $\tw_p =
\rho^\prime(E^{pa}) w_i + \delta_{pa} w_i$. Then using the fact
that $(w_1, \ldots, w_N)$ satisfies (\skews), we obtain
$$ \rho^\prime(E^{cd}) \tw_b + \rho^\prime(E^{bd}) \tw_c$$
$$=\rho^\prime(E^{cd}) \rho^\prime(E^{ba}) w_i + \delta_{ab} \rho^\prime(E^{cd}) w_i
+\rho^\prime(E^{bd}) \rho^\prime(E^{ca}) w_i + \delta_{ac} \rho^\prime(E^{bd}) w_i$$
$$= - \rho^\prime(E^{cd}) \rho^\prime(E^{ia}) w_b + \delta_{ab} \rho^\prime(E^{cd}) w_i
 - \rho^\prime(E^{bd}) \rho^\prime(E^{ia}) w_c + \delta_{ac} \rho^\prime(E^{bd}) w_i$$
$$= - \rho^\prime(E^{ia}) \rho^\prime(E^{cd}) w_b - \delta_{id}  \rho^\prime(E^{ca}) w_b
+ \delta_{ac}  \rho^\prime(E^{id}) w_b + \delta_{ab} \rho^\prime(E^{cd}) w_i$$
$$ - \rho^\prime(E^{ia}) \rho^\prime(E^{bd}) w_c  - \delta_{id}  \rho^\prime(E^{ba}) w_c
+ \delta_{ab}  \rho^\prime(E^{id}) w_c + \delta_{ac} \rho^\prime(E^{bd}) w_i = 0.$$
Thus $(\tw_1, \ldots, \tw_N) \in \PP$. Lemma is now proved.

\end{proof}

\begin{cor} \label{fund}  Let  $(W, \rho)$ be a finite-dimensional
irreducible $\gl_N$-module. If \break
 $L(W,\gamma,h)$ has a critical vector of degree $m \geq 1$ then either $W$  has a fundamental highest weight
$\omega_k$, $1 \leq k \leq N-1$, with respect to $\sln_N$-action, with identity matrix acting with scalar $\alpha = k - mN$, or
$W$ is a 1-dimensional module with identity matrix acting with scalar $\alpha = N - mN$.
\end{cor}

\begin{proof} If the system (\skewtwo) has a non-trivial
solution then the submodule $\tP$ in the tensor module $T(W)$
corresponding to $(W,\rho^\pprime)$ is non-zero. It is a proper
submodule since its component at $q^0$ is trivial. Using the
classification of reducible tensor modules (Theorem \ref{DRc}), we
conclude that $T(W)$ is one of the de Rham modules $\Omega^k
(\T^N)$, $k=1, \ldots, N$. Taking into account the relation
$\alpha = \alpha^\pprime - mN$, we obtain the claim of the
corollary.

\end{proof}

 To complete the proof of Theorem \ref{crith}  it remains to
 establish the following

\begin{lem} \label{hex} If $L(W,\gamma,h)$ has a critical vector that does
not belong the top then $h=0$.

\end{lem}

 \ifnum \value{version}=\value{short} 
\begin{proof}
The statement can be proved using the same method which was used
in the derivation of equation (\Pzb). We will only sketch the idea
of the proof here. A critical vector $g$ is annihilated by $t_0^j
t^r d_0$ for $j\geq 0$. Thus
$$(z^2 d_0 (r,z))_- g = 0.$$
Using (\actfour) we can rewrite that equation explicitly and then project everything to the subspace (\VW), which will
give a system of equations on $f$.
More precisely, one should look at the $r_a$-component in that system.
Since we have by now a fairly detailed description of $f$, the resulting equation can be
simplified yielding $h=0$.

\end{proof}

\else
\begin{proof} Our strategy will be the same as in derivation of
equation (\Pzb). A critical vector $g$ is annihilated by $t_0^j
t^r d_0$ for $j\geq 0$. Thus
$$(z^2 d_0 (r,z))_- g = 0.$$
We will project this equation to the subspace (\VW) in order to derive an equation on $f$.
Finally, we will take $r_a$-component of the resulting equation.
The action of $d_0 (r,z)$
is given by  (\actfour), which has three summands. We will analyze the contribution of each summand
in $z^2 d_0 (r,z)$ separately.

Consider the first summand
$$ - \left( z^2 \left( \sum_{j=1}^\infty \omega_{(-j)} z^{j-1} \right) Y(q^r, z)
+ z^2 Y(q^r, z) \left( \sum_{j=0}^\infty \omega_{(j)} z^{-j-1} \right) \right)_- g . \eqno{(\zera)}$$
We have
$$\omega_{(j)} g = 0 \quad \hbox{\rm for \ } j \geq 2,$$
since $- \omega_{(j)}$ represents $t_0^{j-1} d_0$. Thus the corresponding terms in the above expression may be dropped. We also recall that
$$\omega_{(j)} = \omega^\hyp_{(j)}+\omega^{\gl_N}_{(j)}+\omega^\vir_{(j)} .$$
We further split (\zera) into three summands corresponding to this decomposition. For the case of the Virasoro
field of the hyperbolic lattice component we have
$$\omega^\hyp (z) = \sum_{p=1}^N
\left( \sum_{j=1}^\infty j u_{pj} z^{j-1} \right) \left(
\sum_{k=1}^\infty k v_{pk} z^{k-1} + z^{-1} q_p {\partial \over
\partial q_p} + \sum_{k=1}^\infty  {\partial \over \partial u_{pk}} z^{-k-1}
\right) $$
$$ +  \sum_{p=1}^N
 \left(  \sum_{k=1}^\infty k v_{pk} z^{k-1} + z^{-1} q_p {\partial \over \partial q_p} +  \sum_{k=1}^\infty  {\partial \over \partial u_{pk}} z^{-k-1} \right)
\left(  \sum_{j=1}^\infty  {\partial \over \partial v_{pj}}
z^{-j-1} \right). \eqno{(\zerb)}$$ The first summand in (\zerb)
does not contribute to the projection to (\VW) since it contains
multiplications by $u_{pj}$, while $Y(q^r,z)$ does not involve
differentiations in these variables. Note that we are only
interested in powers $z^j$ in $\omega^\hyp (z)$ with $j\geq -2$.
Thus the only terms that will contribute are:
$$ \sum_{p=1}^N
 \left(  \sum_{k=1}^\infty k v_{pk} z^{k-1}  \right)
\left(  \sum_{j=1}^\infty  {\partial \over \partial v_{pj}}
z^{-j-1} \right).
$$ In operator $Y(q^r,z)$ we may then drop the factors containing
$u_{pj}$ when taking the projection to (\VW). The contribution
that we get will be
$$- \left( z^2 \sum_{p=1}^N
 \left(  \sum_{k=3}^\infty \sum_{j=1}^{k-2}  k v_{pk} {\partial \over \partial v_{pj}}  z^{k-j-2} \right)
\exp \left( - \sum_{p=1}^N r_p \sum_{j=1}^\infty {z^{-j} \over j}
{\partial \over \partial v_{pj}} \right)
 \right)_- f$$
$$- \left( \exp \left( - \sum_{p=1}^N r_p \sum_{j=1}^\infty {z^{-j} \over j} {\partial \over \partial v_{pj}} \right)
\sum_{p=1}^N  \left( z \sum_{j=1}^\infty (j+1) v_{p,j+1} {\partial
\over \partial v_{pj}} +  \sum_{j=1}^\infty j v_{pj}  {\partial
\over \partial v_{pj}}\right) \right)_- f .$$ Let us now take the
$r_a$-coefficient of the expansion in powers of $r$:
$$ \left(  \sum_{p=1}^N
 \left(  \sum_{k=3}^\infty \sum_{j=1}^{k-2}  k v_{pk} {\partial \over \partial v_{pj}}  z^{k-j} \right)
\left(  \sum_{j=1}^\infty {z^{-j} \over j} {\partial \over
\partial v_{aj}} \right) f \right)_- $$
$$+ \left( \left(  \sum_{j=1}^\infty {z^{-j} \over j} {\partial \over \partial v_{aj}} \right)
\sum_{p=1}^N  \left( z \sum_{j=1}^\infty (j+1) v_{p,j+1}
{\partial \over
\partial v_{pj}} +  \sum_{j=1}^\infty j v_{pj}  {\partial \over
\partial
v_{pj}} \right)  f  \right)_- .  \eqno{(\zerd)}$$ Since $f$ is
linear in $v_{pm}$,
$$f = \sum_{p=1}^N v_{pm} \otimes w_p \otimes \vh, \quad w_p \in W,$$
the first summand in (\zerd) vanishes, while the second simplifies to
$$ \left( z^{-m} \sum_{j=1}^\infty   {\partial \over \partial v_{aj}}
+  z^{-m} \sum_{j=1}^\infty  {\partial \over \partial v_{aj}}
\right) f
$$
$$ = 2 z^{-m} 1 \otimes w_a \otimes \vh . \eqno{(\zere)}$$

Next, let us consider the contribution of the Virasoro field of $\Vgl$:
$$ - \left( z^2 \left( \sum_{j=1}^\infty \omega^{\sgl_N}_{(-j)} z^{j-1} \right) Y(q^r, z)
+ z^2 Y(q^r, z) \left( z^{-1} \omega^{\sgl_N}_{(0)} + z^{-2} \omega^{\sgl_N}_{(1)}\right) \right)_- g .
\eqno{(\zerf)}$$
The operators $\omega^{\sgl_N}_{(j)}$ with $j \leq 0$ increase the degree in the component $\Lgl$
and thus do not contribute to the projection to the space (\VW), and only the term with $\omega^{\sgl_N}_{(1)}$
will contribute.  The Virasoro field of $\hgl_N$ is a sum of the Virasoro fields of $\hsl_N$ (\virfa)  and the Heisenberg algebra (\virfb).

Using (\virfa) we can write
$$\omega^{\ssl_N}_{(1)} =
{1 \over 2(N+1)} \left( \sum_{k=1}^\infty \sum_{i,j=1}^N E^{ij}_{(-k)} E^{ji}_{(k)}
+ \sum_{k=0}^\infty \sum_{i,j=1}^N E^{ji}_{(-k)} E^{ij}_{(k)} \right.$$
$$ \left. - {1 \over N} \sum_{k=1}^\infty I_{(-k)} I_{(k)}
- {1 \over N}  \sum_{k=0}^\infty  I_{(-k)} I_{(k)} \right), $$
and the terms that contribute to the projection are
$${1 \over 2(N+1)} \left(  \sum_{i,j=1}^N E^{ij}_{(0)} E^{ji}_{(0)} - {1 \over N}  I_{(0)} I_{(0)} \right),$$
which is a multiple of the Casimir operator for $\sln_N$. If $W$ corresponds to the tensor module of
$k$-forms, $k=1,\ldots,N$, this operator will act on the space (\VW) with scalar
${k (N-k) \over 2 N}$ (see (\casim), this also includes the case of a trivial $\sln_N$ module $W$ when $k=N$).

 Analogously, for the Virasoro field (\virfb) of the Heisenberg algebra $\omega^\shei_{(1)}$ will
be acting on $f$ with the scalar
$${1 \over 2N}  I_{(0)} I_{(0)} - {1\over 2} I_{(0)} = {(k-mN)^2 \over 2N} - {k-mN \over 2}.$$

Going back to (\zerf), we get the contribution of the Virasoro field in $\Vgl$.
$$ -\left( {k (N-k) \over 2 N} + {(k-mN)^2 \over 2N} - {k-mN \over 2} \right)
 \exp \left( - \sum_{p=1}^N r_p \sum_{j=1}^\infty {z^{-j} \over j} {\partial \over \partial v_{pj}} \right) f,$$
and its $r_a$-term will yield
$${z^{-m} \over m} \left( {k (N-k) \over 2 N} + {(k-mN)^2 \over 2N} - {k-mN \over 2} \right) 1 \otimes w_a \otimes \vh . \eqno{(\zerk)}$$

Now let us deal with the Virasoro field of $\Lvir$. The corresponding term is
$$ - \left( z^2 \left( \sum_{j=1}^\infty \omega^\vir_{(-j)} z^{j-1} \right) Y(q^r, z)
+ z^2 Y(q^r, z) \left( z^{-1} \omega^\vir_{(0)} + z^{-2} \omega^\vir_{(1)}\right) \right)_- g .
\eqno{(\zerh)}$$
Since the operators $\omega^\vir_{(j)}$ with $j \leq 0$ increase the degree in the component $\Lvir$,
the only term that contributes to the projection to (\VW) is $\omega^\vir_{(1)}$, which acts on (\VW)
with scalar $h$. Thus the $r_a$-term of (\zerh) gives the contribution
$$ \left( \sum_{j=1}^\infty {z^{-j} \over j} {\partial \over \partial v_{aj}} \right) \omega^\vir_{(1)} f =
{z^{-m} \over m} h \; 1 \otimes w_a \otimes \vh . \eqno{(\zerl)}$$

Next we shall look at the the summand $-\sum\limits_{a,b=1}^N r_a u^b (z) E^{ab}(z) Y(q^r,z)$ in (\actfour). Its $r_a$-term
is
$$-\sum\limits_{p=1}^N  \left(z^2  \left( \sum_{j=1}^\infty j u_{pj} z^{j-1}
+  \sum_{j=1}^\infty  {\partial \over \partial v_{pj}} z^{-j-1}
\right) \times \left( \sum_{k \in \Z} E^{ap}_{(k)} z^{-k-1}
\right) \right)_- g .$$ When we consider the projection to (\VW),
we can drop terms with multiplications by $u_{pj}$ and
$E^{ap}_{(k)}$ with $k\leq -1$, while for  $E^{ap}_{(k)}$ with $k
\geq 1$ we may use (\Eqq), which yields
$$ -\sum\limits_{p=1}^N  \left(  \sum_{j=1}^\infty  {\partial \over \partial v_{pj}} z^{-j}  \right) E^{ap}_{(0)} f$$
$$ -\sum\limits_{p=1}^N  \left(  \sum_{j=1}^\infty  {\partial \over \partial v_{pj}} z^{-j}  \right)
 \left( \sum_{k=1}^\infty \left( \sum_{i=1}^{k-1} i v_{pi} z^i + q_p {\partial \over \partial q_p} \right)
{z^{-k} \over k} {\partial \over \partial v_{ak}} \right) f .$$
Taking into account that $f$ is linear in $v_{pm}$, this can be
simplified to the following
$$ -\sum\limits_{p=1}^N  \left(  \sum_{j=1}^\infty  {\partial \over \partial v_{pj}} z^{-j}  \right) E^{ap}_{(0)} f
 -  N \left( \sum_{k=1}^\infty \left( \sum_{i=1}^{k-1} i  \right)
{z^{-k} \over k} {\partial \over \partial v_{ak}} \right) f $$
$$ =  - z^{-m} \sum\limits_{p=1}^N 1 \otimes E^{ap} w_p \otimes \vh - z^{-m} N {m-1 \over 2} 1 \otimes
w_a \otimes \vh . \eqno{(\zeri)}$$
Lemma \ref{redten} provides a relation between components $(w_1, \ldots, w_N)$ and submodules in the
tensor modules $\Omega^k (\T^N)$. Using computations in the tensor module of $k$-forms one can
show that
$$\sum\limits_{p=1}^N  E^{ap} w_p = (N - k - m + 1) w_a ,$$
thus (\zeri) reduces to
$$  - z^{-m} \left(  (N - k - m + 1) + N {m-1 \over 2} \right)  1 \otimes w_a \otimes \vh . \eqno{(\zerj)}$$

The $r_a$-coefficient of the last summand
$$ \sum_{p=1}^N \left( z^2 \left( {d \over dz} u^p (z) \right) Y(q^r,z) \right)_- g$$
is
$$\left( \sum_{j=1}^\infty j (j-1) u_{aj} z^{j} + \sum_{j=1}^\infty (-j-1)  {\partial \over \partial v_{aj}} z^{-j}
 \right)_- g, $$
and its projection to (\VW) yields
$$ \sum_{j=1}^\infty (-j-1)  {\partial \over \partial v_{aj}} z^{-j} f = - (m+1) z^{-m} 1 \otimes w_a \otimes \vh
\eqno{(\zerm)} $$
Finally, collecting (\zere),  (\zerk), (\zerl), (\zerj) and (\zerm) together, we get
$$ 0 = \left( 2  + {1 \over m} \left( {k (N-k) \over 2 N} + {(k-mN)^2 \over 2N} - {k-mN \over 2} \right)
+ {h \over m} -  (N - k - m + 1) \right.$$
$$\hbox{\hskip 2cm} \left.  - N {m-1 \over 2} - (m+1) \right) w_a
 = {h \over m}  w_a .$$
Since $a$ is arbitrary, we can choose it so that $w_a \neq 0$.
Thus $h=0$, which was to be demonstrated.

\end{proof}
\fi

\

\section{Chiral de Rham complex}

Chiral de Rham complex was introduced by Malikov et al. in
\cite{MSV}. In case a of torus $\T^N$ the space of this
differential complex is a tensor product of two vertex (super)
algebras
$$\Vhyp \otimes \VZ. $$
Here $\VZ$ is the lattice vertex superalgebra of the standard euclidean lattice $\Z^N$. Before we define the
differential of this complex, let us review the structure of $\VZ$. The vertex superalgebra $\VZ$ has two
main realizations -- the bosonic realization and the fermionic one, with boson-fermion correspondence being
an isomorphism between the two models. For our purposes it will be more convenient to use the fermionic
realization of $\VZ$.

Consider the Clifford Lie superalgebra $\Cl_N$ of ``charged free fermions'' with basis

$$\{ \varphi^p_{(j)},  \psi^p_{(j)} |{p=1,\ldots,N},{j\in\Z} \}$$ of
its odd part and a 1-dimensional even part spanned by a central
element $K$. The Lie bracket in $\Cl_N$ is given by
$$[ \varphi^a_{(m)}, \psi^b_{(n)} ] = \delta_{ab} \delta_{m, -n-1} K, \quad
[ \varphi^a_{(m)}, \varphi^b_{(n)} ] = [ \psi^a_{(m)}, \psi^b_{(n)} ] = 0 .$$
Define formal fields
$$\varphi^a (z) = \sum_{j\in \Z} \varphi^a_{(j)} z^{-j-1} , \quad\quad
\psi^a (z) = \sum_{j\in \Z} \psi^a_{(j)} z^{-j-1} , \quad\quad
K(z) = K z^0 .$$
With this choice of fields $\Cl_N$ becomes a vertex Lie superalgebra since the only non-trivial relation
between these fields is
$$ \left[ \varphi^a (z_1), \psi^b (z_2) \right] = \delta_{ab} K(z_2) \left[ z_1^{-1} \delta \left( {z_2 \over z_1}
\right) \right] .$$
The lattice vertex superalgebra $\VZ$ is isomorphic to the universal enveloping vertex algebra of $\Cl_N$ at
level $1$. As a vector space it is the unique $\Cl_N$-module generated by vacuum vector $\vac$, satisfying
$$K \vac = \vac, \quad \varphi^p_{(j)} \vac =  \psi^p_{(j)} \vac = 0 \hbox{{\rm \ for \ }} j \geq 0, \;
p = 1, \ldots, N.$$
In its fermionic realization $\VZ$ is the exterior algebra with generators
$\{ \varphi^p_{(j)},  \psi^p_{(j)} |^{p=1,\ldots,N}_{j \leq -1} \}$ and is irreducible as a module over $\Cl_N$.
The state-field correspondence map $Y$ is given by the standard formula (\vla).

We fix the Virasoro element in $\VZ$:
$$\om^\fer = \sum_{p=1}^N \varphi^p_{(-2)} \psi^p_{(-1)} \vac .$$
The rank of this VOA is $-2N$.

 It is well-known that vertex superalgebra $\VZ$ contains a level $1$ simple $\hgl_N$ vertex algebra.
The fields generating this subalgebra are
$$E^{ab} (z) = : \varphi^a (z) \psi^b (z) : .$$
It is easy to check that these satisfy relations (\glbrak) and the central element of  $\hgl_N$ acts as identity
operator. It is also straightforward to verify that the Virasoro element (\glvir) in the $\hgl_N$ vertex algebra
maps to $\om^\fer$ under this embedding.

 Let us define two $\Z$-gradings on $\VZ$. The {\it fermionic} degree is defined by
$$\deg_\fer (\varphi^p_{(j)}) = 1, \quad \deg_\fer (\psi^p_{(j)}) = -1, \quad
\deg_\fer (K) = \deg_\fer (\vac) = 0 .$$

The {\it bosonic} grading is defined as follows:
$$\deg_\bos (\varphi^p_{(j)}) = -j-1, \quad \deg_\bos (\psi^p_{(j)}) = -j, \quad
\deg_\bos (K) = \deg_\bos (\vac) = 0 .$$

Let $\VZ^k$ be the subspace of the elements of fermionic degree $k$. We have a decomposition
$$\VZ = \mathop\oplus\limits_{k \in \Z} \VZ^k .$$
Note that each subspace $\VZ^k$ is a $\hgl_N$-submodule, which is
graded by the bosonic degree. Its structure is described by the
following well-known result (see e.g. \cite{Fr} or \cite{KL}):

\begin{thm} \label{wellknown} For each $k\in \Z$, $\VZ^k$ is an irreducible
$\hgl_N$-module at level $1$. Let ${\widetilde V}_{\Z^N}^k$ be the
non-trivial component of $\VZ^k$ of the lowest bosonic degree. If
$k = 0 \; \mod \; N$ then ${\widetilde V}_{\Z^N}^k$ is
$1$-dimensional. If $k = k^\prime \; \mod \; N$ with $1 \leq
k^\prime < N$, then as an $\sln_N$-module ${\widetilde
V}_{\Z^N}^k$ has the fundamental highest weight $\om_{k^\prime}$.
The identity matrix of $\gl_N$ acts on ${\widetilde V}_{\Z^N}^k$
as $k \, \Id$.
\end{thm}

 Combining this result with Theorem \ref{thDK}, we get

\begin{cor} \label{chmd} The space
$$\Mhyp (\gamma) \otimes \VZ^k$$
has a structure of a module for the Lie algebra $\Venone$ of
vector fields.

\end{cor}

For these modules the Virasoro tensor factor $\Lvir (h)$ is $1$-dimensional ($h=0$). The modules in this family
are precisely the exceptional modules $L(W, \gamma,  h)$ for which Theorem \ref{irr} \ does not
claim irreducibility. We are going to see below that these modules are in fact reducible.

 Let us express the action (\actthree), (\actfour) of the Lie algebra $\Venone$ on
$\Mhyp (\gamma) \otimes \VZ^k$ using the fermionic realization:
$$d_a(r,z) \mapsto Y(d_a(r), z), \quad d_0(r,z) \mapsto Y(d_0(r), z), $$
where
$$d_a (r) = v^a_{(-1)} q^r + \sum_{p=1}^N r_p \varphi^p_{(-1)} \psi^a_{(-1)} q^r ,$$
$$d_0 (r) = -  \left( \om^\hyp_{(-1)} q^r + \om^\fer_{(-1)} q^r
+ \sum_{a,b=1}^N r_a u^b_{(-1)}  \varphi^a_{(-1)} \psi^b_{(-1)} q^r
- \sum_{p=1}^N r_p u^p_{(-2)} q^r  \right) $$
$$=  - \left( \sum_{p=1}^N u^p_{(-1)} v^p_{(-1)} q^r
+ \sum_{p=1}^N \varphi^p_{(-2)} \psi^p_{(-1)} q^r
+ \sum_{a,b=1}^N r_a u^b_{(-1)}  \varphi^a_{(-1)} \psi^b_{(-1)} q^r
\right) .$$
Here we used the relation
$$ \om^\hyp_{(-1)} q^r =  \sum_{p=1}^N u^p_{(-1)} v^p_{(-1)} q^r
+ \sum_{p=1}^N r_p u^p_{(-2)} q^r .$$

 Following \cite{MSV}, let us now introduce the differential
$$ \ldots {\darrow} \Mhyp (\gamma) \otimes \VZ^k
{\darrow} \Mhyp (\gamma) \otimes \VZ^{k+1}
{\darrow} \ldots $$
of the chiral de Rham complex.

Let
$$Q = \sum\limits_{p=1}^N v^p_{(-1)} \varphi^p_{(-1)} \vac$$
and set $\dd = Q_{(0)}$, i.e., $\dd$ is a coefficient at $z^{-1}$ in
$Y(Q,z) =  \sum\limits_{p=1}^N v^p (z) \varphi^p (z)$. Vanishing
of the supercommutator
$$ \left[ Y(Q, z_1), Y(Q, z_2) \right] = 0$$
implies $\dd \circ \dd = 0$.

\begin{thm} \label{chhom} The map
$$\dd: \;  \Mhyp (\gamma) \otimes \VZ^k
\rightarrow \Mhyp (\gamma) \otimes \VZ^{k+1} $$ is a homomorphism
of $\Venone$-modules.
\end{thm}

The statement of this theorem is equivalent to the claim that the
following operators on $ \Mhyp (\gamma) \otimes \VZ$ commute:
$$\left[ \dd,  d_a (r,z) \right] = 0,$$
$$\left[ \dd,  d_0 (r,z) \right] = 0.$$
The proof of these relations will be based on the following simple observation:

\begin{lem} \label{compcom} Let $V$ be a vertex superalgebra and let
$a,b \in V$. Suppose that $a_{(0)} b = 0$. Then
$$\left[ a_{(0)} ,  Y(b,z) \right] = 0.$$
\end{lem}

\begin{proof} Since  $a_{(0)} b = 0$, the commutator
formula (\comm) yields
$$\left[ Y(a,z_1) ,  Y(b,z_2) \right] =
\sum_{j \geq 1}  {1 \over j!} Y(a_{(j)} b, z_2) \left[ z_1^{-1}
\left( {\partial \over \partial z_2} \right)^j \delta \left( {z_2
\over z_1} \right) \right] .$$ However the right hand side does
not contain terms with $z_1^{-1}$ and the claim of the lemma
follows.

\end{proof}

 Let us continue with the proof of the theorem. We need to show that
$$Q_{(0)} d_a (r) = 0 \quad \hbox{\rm and} \quad Q_{(0)} d_0 (r) = 0.$$
Since $Y(Q,z) = \sum\limits_{i=1}^N :v^i (z) \varphi^i (z):$, we have
$$Q_{(0)} = \sum_{i=1}^N \left(
\sum_{j=0}^\infty \varphi^i_{(-j-1)} v^i_{(j)} + \sum_{j=1}^\infty v^i_{(-j)} \varphi^i_{(j-1)}
\right).$$
It is easy to see that $v^i_{(j)} d_a(r) = 0$ for $j \geq 1$ and
$\varphi^i_{(j)}  d_a(r) = 0$ for $j \geq 1$. Thus
$$Q_{(0)}  d_a(r) =  \sum_{i=1}^N \left(\varphi^i_{(-1)} v^i_{(0)} +
 v^i_{(-1)} \varphi^i_{(0)} \right) d_a(r)$$
$$=  \sum_{i=1}^N r_i \varphi^i_{(-1)} v^a_{(-1)} q^r
- \sum_{i=1}^N \sum_{p=1}^N r_p  v^i_{(-1)} \varphi^p_{(-1)} \varphi^i_{(0)}
\psi^a_{(-1)} q^r$$
$$=  \sum_{i=1}^N r_i \varphi^i_{(-1)} v^a_{(-1)} q^r
-  \sum_{p=1}^N r_p v^a_{(-1)} \varphi^p_{(-1)}  q^r = 0.$$

Let us now show that $Q_{(0)} d_0 (r) = 0$.
Since  $v^i_{(j)} d_0(r) = 0$ for $j \geq 2$ and
$\varphi^i_{(j)}  d_0(r) = 0$ for $j \geq 1$, we get
$$- Q_{(0)} d_0 (r) = \sum_{i=1}^N \left( \varphi^i_{(-1)} v^i_{(0)} +
\varphi^i_{(-2)} v^i_{(1)} + v^i_{(-1)} \varphi^i_{(0)} \right) (-d_0(r)).$$
Let us compute each of three terms in the right hand side separately:
$$\left(  \sum_{i=1}^N  \varphi^i_{(-1)} v^i_{(0)}  \right)  (-d_0(r)) $$
$$=  \sum_{i=1}^N  \sum_{p=1}^N  \varphi^i_{(-1)}  v^i_{(0)} u^p_{(-1)} v^p_{(-1)} q^r
+  \sum_{i=1}^N  \sum_{p=1}^N  \varphi^i_{(-1)}  v^i_{(0)} \varphi^p_{(-2)} \psi^p_{(-1)} q^r$$
$$+  \sum_{i=1}^N  \sum_{a,b=1}^N
r_a \varphi^i_{(-1)}  v^i_{(0)} u^b_{(-1)} \varphi^a_{(-1)}  \psi^b_{(-1)} q^r$$
$$=  \sum_{i=1}^N  \sum_{p=1}^N  r_i \varphi^i_{(-1)} u^p_{(-1)} v^p_{(-1)} q^r
+  \sum_{i=1}^N  \sum_{p=1}^N  r_i \varphi^i_{(-1)} \varphi^p_{(-2)} \psi^p_{(-1)} q^r$$
$$+  \sum_{i=1}^N  \sum_{a=1}^N \sum_{b=1}^N
r_a r_i u^b_{(-1)} \varphi^i_{(-1)}  \varphi^a_{(-1)}  \psi^b_{(-1)} q^r. \eqno{(\chone)} $$
The last summand in (\chone) vanishes since it is antisymmetric in $\{ a, i \}$. Next,
$$\left(  \sum_{i=1}^N  \varphi^i_{(-2)} v^i_{(1)}  \right)  (-d_0(r)) $$
$$= \sum_{i=1}^N  \sum_{p=1}^N \varphi^i_{(-2)} v^i_{(1)} u^p_{(-1)} v^p_{(-1)} q^r
+  \sum_{i=1}^N  \sum_{p=1}^N  \varphi^i_{(-2)}  v^i_{(1)} \varphi^p_{(-2)} \psi^p_{(-1)} q^r$$
$$+  \sum_{i=1}^N  \sum_{a,b=1}^N
r_a \varphi^i_{(-2)}  v^i_{(1)} u^b_{(-1)} \varphi^a_{(-1)}  \psi^b_{(-1)} q^r$$
$$=  \sum_{p=1}^N \varphi^p_{(-2)}  v^p_{(-1)} q^r
+  \sum_{i=1}^N  \sum_{a=1}^N
r_a \varphi^i_{(-2)}  \varphi^a_{(-1)}  \psi^i_{(-1)} q^r . \eqno{(\chtwo)}$$
And finally,
$$\left(  \sum_{i=1}^N  v^i_{(-1)} \varphi^i_{(0)}  \right)  (-d_0(r)) $$
$$= \sum_{i=1}^N  \sum_{p=1}^N  v^i_{(-1)} \varphi^i_{(0)} u^p_{(-1)} v^p_{(-1)} q^r
+  \sum_{i=1}^N  \sum_{p=1}^N   v^i_{(-1)} \varphi^i_{(0)} \varphi^p_{(-2)} \psi^p_{(-1)} q^r$$
$$+  \sum_{i=1}^N  \sum_{a,b=1}^N
r_a v^i_{(-1)} \varphi^i_{(0)} u^b_{(-1)} \varphi^a_{(-1)}  \psi^b_{(-1)} q^r$$
$$ = - \sum_{i=1}^N  \sum_{p=1}^N   v^i_{(-1)} \varphi^p_{(-2)}  \varphi^i_{(0)} \psi^p_{(-1)} q^r
- \sum_{i=1}^N  \sum_{a,b=1}^N
r_a v^i_{(-1)} u^b_{(-1)} \varphi^a_{(-1)}  \varphi^i_{(0)} \psi^b_{(-1)} q^r$$
$$= - \sum_{p=1}^N   v^p_{(-1)} \varphi^p_{(-2)}  q^r
- \sum_{i=1}^N  \sum_{a=1}^N r_a v^i_{(-1)} u^i_{(-1)}
\varphi^a_{(-1)} q^r .\eqno{(\chthree)}$$ Combining (\chone),
(\chtwo) and (\chthree) we get $Q_{(0)} (-d_0 (r)) = 0$, and the
theorem is proved.

\

 Let us present here a diagram of the Chiral de Rham complex for $N=2$. On the diagram, the fermionic degree increases in the horizontal direction and bosonic in vertical.

\

\hskip 25pt
\vbox{
\hbox{$\hskip 117pt \Omega^0 \hskip 26pt \Omega^1 \hskip 26pt \Omega^2$}
\vskip -4pt
\hbox{$\hskip 115pt \scriptstyle k=0 \hskip 22pt k=1 \hskip 23pt k=2$}
\vskip -4pt
\hbox{$\hskip 120pt \bullet \marrow \bullet  \marrow \bullet \rightarrow$}
\vskip -4pt
\hbox{$\hskip 121pt \vrt \hskip 36.5pt \vrt \hskip 36.5pt \vrt$}
\vskip -15pt
\hbox{$\hskip 73pt \scriptstyle k=-1 \hskip 5pt \dd \hskip 105pt \dd \hskip 10pt \scriptstyle k=3$}
\vskip -5pt
\hbox{$\hskip 82pt \bullet \marrow \bullet  \hskip 30.5pt \bullet \hskip 30.5pt \bullet \marrow \bullet \rightarrow$}
\vskip -3pt
\hbox{$\hskip 84pt \vrt \hskip 36pt \vrt \hskip 36.5pt  \vrt \hskip 36.5pt \vrt \hskip 36.5pt \vrt$}
\vskip -15pt
\hbox{$\hskip 37pt \scriptstyle k=-2 \hskip 5pt \dd \hskip 178pt \dd \hskip 10pt \scriptstyle k=4$}
\vskip -5pt
\hbox{$\hskip 45pt \bullet \marrow \bullet  \hskip 30pt \bullet \hskip 30pt \bullet  \hskip 30.5pt \bullet \hskip 30.5pt  \bullet \marrow \bullet \rightarrow$}
\vskip -3pt
\hbox{$\hskip 47pt \vrt \hskip 36pt \vrt \hskip 36pt  \svrt \hskip 36.5pt  \svrt \hskip 36.5pt \svrt \hskip 37pt \vrt \hskip 36pt \vrt$}
\vskip -8pt
\hbox{$\hskip 45 pt \bullet \hskip 30pt \bullet  \hskip 143pt \bullet \hskip 32pt \bullet   \rightarrow$}
\vskip -2pt
\hbox{$\hskip 47pt \vrt \hskip 36pt \svrt \hskip 149pt \svrt \hskip 36pt \vrt$}
\vskip -15pt
\hbox{$\scriptstyle k=-3 \hskip 5pt \dd \hskip 252pt \dd \hskip 8pt \scriptstyle k=5$}
\vskip -5pt
\hbox{$\hskip 8pt \bullet \marrow \bullet  \hskip 217pt  \bullet \marrow \bullet \rightarrow$}
\vskip -4pt
\hbox{\hskip 9pt \svrt \hskip 37pt \svrt \hskip 223pt \svrt \hskip 34pt \svrt}
}

\

The tops of the modules $\Mhyp(\gamma) \otimes \VZ^k$ with $0 \leq
k \leq N$ are the spaces $q^\gamma \Omega^k(\T^N)$ of differential
$k$-forms that form the classical de Rham complex. Non-trivial
$\VenN$-submodules in these tops generate non-trivial $\Venone$
submodules in corresponding modules $\Mhyp(\gamma) \otimes \VZ^k$.

It was proved in \cite{MSV} that the cohomology of the chiral de
Rham complex coincides with the classical de Rham cohomology. This
implies, in particular, that for $k<0$ or $k >N$ the short
sequences
$$\Mhyp(\gamma) \otimes \VZ^{k-1} \darrow
\Mhyp(\gamma) \otimes \VZ^k \darrow \Mhyp(\gamma) \otimes \VZ^{k+1}$$
are exact. Using this fact, we get

\begin{cor} \label{CG} (i) For $k \leq 0$, $\Venone$-modules $\Mhyp(\gamma)
\otimes \VZ^k$ have non-trivial critical vectors.

(ii) For  $k \geq N$, $\Venone$-modules $\Mhyp(\gamma) \otimes
\VZ^k$ are not generated by their top spaces.
\end{cor}

\begin{proof} We can see from the above diagram that for $k<0$
the images of the top vectors in $\Mhyp(\gamma) \otimes \VZ^k$ are
non-trivial critical vectors in $ \Mhyp(\gamma) \otimes
\VZ^{k+1}$. For $k \geq N$, the top spaces of $\Mhyp(\gamma)
\otimes \VZ^k$ are in the kernel of $\dd$. Thus the submodules
generated by the tops are annihilated by $\dd$ as well. Since the
map $\dd$ is non-zero, these submodules are proper.
\end{proof}

As a result we see that all modules that belong to the chiral de
Rham complex are reducible. The claim of Corollary \ref{CG} is
consistent with the existence of the contragredient pairing given
by Theorem \ref{duality}:
$$ \left( \Mhyp(\gamma) \otimes \VZ^k \right) \times
 \left( \Mhyp(\gamma) \otimes \VZ^{N-k} \right) \rightarrow \C .$$
For the chiral de Rham complex this duality was constructed in
\cite{MS2}.

\section{Acknowledgements}
The first author is supported in part by a grant from the Natural Sciences
and Engineering Research Council of Canada.  
The second author is supported in part by the CNPq grant
(301743/2007-0) and by the Fapesp grant (2010/50347-9).
Part of this work was carried out during the visit of the first author
to the University of S\~ao Paulo in 2009. This author would like
to thank the University of S\~ao Paulo for hospitality and
excellent working conditions and Fapesp (2008/10471-2) for
financial support.

\end{document}